\documentclass[12pt,twoside,reqno]{amsart}
\usepackage{amsxtra,amscd}
\usepackage{graphicx}
\usepackage{amsmath}
\usepackage{amsfonts}
\usepackage{amssymb}

\newtheorem{theorem}{Theorem}[section]
\newtheorem{lemma}[theorem]{Lemma}
\newtheorem{maintheorem}[theorem]{Main Theorem}
\theoremstyle{definition}
\newtheorem{definition}[theorem]{Definition}

\newtheorem{proposition}[theorem]{Proposition}
\newtheorem{remark}[theorem]{Remark}

\newtheorem{corollary}[theorem]{Corollary}

\numberwithin{equation}{section}

\begin{document}
\title{the automorphism groups of
quasi-galois closed arithmetic schemes}
\author{Feng-Wen An}
\address{School of Mathematics and Statistics, Wuhan University, Wuhan,
Hubei 430072, People's Republic of China} \email{fwan@amss.ac.cn}

\begin{abstract}
Assume that $X$ and $Y$ are arithmetic schemes, i.e., integral
schemes of finite types over $Spec(\mathbb{Z})$. Then $X$ is said to
be quasi-galois closed over $Y$ if $X$ has a unique conjugate over
$Y$ in some certain algebraically closed field, where the conjugate
of $X$ over $Y$ is defined in an evident manner. Now suppose that
$\phi:X\rightarrow Y$ is a surjective morphism of finite type such
that  $X$ is quasi-galois closed over $Y$. In this paper the main
theorem says that the function field $ k\left( X\right)$ is
canonically a Galois extension of $k\left( Y\right)$ and the
automorphism group ${Aut}(X/Y)$ is isomorphic to the Galois group
$Gal(k(X)/k\left( Y\right) )$; in particular, $\phi$ must be affine. Moreover, let $\dim X=\dim Y$.
Then $X$ is a pseudo-galois cover of  $Y$ in the sense of Suslin-Voevodsky.
\end{abstract}

\maketitle





{\tiny{
\begin{center}
{Contents}
\end{center}

\qquad {Introduction}

\qquad {1. Notation and Definition}

\qquad\quad {1.1. Convention}

\qquad\quad {1.2. Affine covering with values in a given field}

\qquad\quad {1.3. Quasi-galois closed varieties}

\qquad {2. Statement of The Main Theorem}

\qquad {3. Proof of The Main Theorem}

\qquad\quad {3.1. Affine structures}

\qquad\quad {3.2. A quasi-galois closed variety has only one maximal affine structure among others

\qquad\quad\qquad with values in a fixed field}

 \qquad\quad {3.3. Definition for conjugations of a given field}

 \qquad\quad {3.4. A quasi-galois field has only one conjugation}

\qquad\quad {3.5. Definition for conjugations of an open set}

\qquad\quad {3.6. Quasi-galois closed varieties and conjugations of open sets}

\qquad \quad {3.7. Automorphism groups of quasi-galois closed varieties and Galois groups of the

\qquad\quad\qquad function fields}

\qquad \quad {3.8. Proof of the main theorem}

\qquad {References}}}

\section*{Introduction}

Let $k_{1}$ be an algebraic extension of a field $k$ and let $X$ be
an algebraic variety defined over $k_{1}$. Then by
a $k-$automorphism $\sigma$ of $\overline{k}$, we get a conjugate
$X^{\sigma}$ of $X$ over $k$. $X$ is said to be
\emph{normally algebraic} over $k$, defined by Weil, if $X$ coincides with each of the
conjugates of $X$ over $k$ (see \cite{weil}).
It is well-known that algebraic varieties and their conjugates  behave like conjugates of
fields and have almost all of the algebraic properties (for example, see \cite{SGA1,weil}). But their complex
topological properties are very different from each other (for example, see \cite{Milne,Serre}).
On the other hand, in the geometric version of class field theory,  following Weil's
\cite{weil2} and \cite{weil3}, Lang uses algebraic varieties to describe unramified class
fields over function fields in several variables (see \cite{Lang}); then in virtue of Bloch's \cite{Bloch}  and
others' foundations, Kato and Saito use algebraic fundamental groups
to obtain unramified class field theory (see \cite{Kato,Saito}). Here, the main feature is to use the abelianized fundamental groups of
algebraic or arithmetic schemes to describe
abelian class fields (for example, see \cite{Schmidt,Raskind,w1,w2}).

Motivated by those works, in this paper we will suggest a definition that
 an arithmetic variety is said to be
\emph{quasi-galois closed} if it has a unique conjugate in an algebraically
closed field (see \emph{Definition 1.1}), which can be regarded as a generalization from
 the notion that an algebraic variety is normally
algebraic over a number field to the one that an arithmetic variety is
quasi-galois closed over a fixed arithmetic one.  Here, an
\emph{arithmetic variety} is an integral scheme of finite type over $
Spec\left( \mathbb{Z}\right) $. Then we will try to use these
relevant data of such arithmetic varieties to obtain some information of Galois
extensions of function fields in several variables.

The following is the \emph{Main Theorem} of the paper (i.e., \emph{Theorem 2.1}).

\begin{maintheorem}
\textbf{\emph{(Theorem 2.1)}} Let $X$ and $Y$ be two arithmetic
varieties. Assume that $X$ is quasi-galois closed over $Y$ by a
surjective morphism $\phi$ of finite type. Then there are the following statements.

\begin{itemize}
\item $f$ is affine.

\item $k\left( X\right) $ is canonically a Galois
extension of $k(Y)$.

\item There is a group isomorphism
$$
{Aut}\left( X/Y\right) \cong Gal(k\left( X\right) /k(Y)).
$$

\item Particularly, let $\dim X=\dim Y$.
Then  $X$ is a pseudo-galois cover of  $Y$ in the sense of Suslin-Voevodsky.
\end{itemize}
\end{maintheorem}

See \cite{VS1} for the definition of \emph{pseudo-galois} covers of schemes. Note that here $k(X)$ is not necessarily algebraic
over $k(Y)$ by $\phi$ in the first property above. That is, the morphism $\phi$
is not necessarily finite.

Hence, the \emph{Main Theorem} of the paper shows us some evidence that  there exists a
nice relationship between quasi-galois closed arithmetic varieties
and Galois extensions of functions fields in several variables.

For the case that $\phi$ is finite, it can be seen that quasi-galois closed arithmetic varieties behave like
Galois extensions of number fields and their automorphism groups can
be regarded as the Galois groups of the field extensions.

In deed, one has been
attempted to use the data of such varieties $X/Y$ to describe a
given Galois extension $E/F$ for a long time and one says that $X/Y$ are a \emph{model} for
$E/F$ if the Galois group $Gal\left( E/F\right) $ is isomorphic to
the automorphism group ${Aut}\left( X/Y\right)$ (for example, see \cite{SGA1,GIT,Raskind,VS1,SV2}).
The \emph{Main Theorem} gives us such a model for function fields in several variables.

In \cite{VS1,SV2}, Suslin and Voevodsky obtain several good properties for pseudo-galois covers
of varieties for the case that the morphism
$\phi$ is finite, where they also give the existence of pseudo-galois covers. If the arithmetic varieties are of the same dimensions, it is seen that there is no essential difference between our \textquotedblleft quasi-galois closed\textquotedblright\ and \textquotedblleft pseudo-galois cover\textquotedblright.

However, there is a main difference between the two types of covers if the structure morphism is not finite. For example, let $t$ be a variable over $\mathbb{Q}$. Then $Spec(\mathbb{Z}[t])/Spec(\mathbb{Z})$ is quasi-galois closed but not pseudo-galois.
Hence, to some degree, the \emph{Main Theorem} of the paper gives us a sufficient
condition for the existence of such a pseudo-galois cover in a more
generalized case in the category of arithmetic varieties, where the function fields are in several variables.

The \emph{Main Theorem} of the paper can be regarded as a generalization of
\emph{Proposition 1.1} in \cite{SGA1}, \emph{Page 106}, for the case of function fields in several variables.

Now let us give some applications of quasi-galois closed covers such as the following.

In \cite{An3} we will prove the existence of quasi-galois closed covers of arithmetic schemes and then by these covers we will give an explicit construction of the geometric model for a prescribed Galois extension of a function field in several variables over a number field.

In \cite{An5} we will use quasi-galois closed covers to define and compute a qc fundamental group for an arithmetic scheme. Then we will prove that the \'{e}tale fundamental group of an arithmetic scheme is a normal subgroup in our qc fundamental group. Hence, our group gives us a prior estimate of the \'{e}tale fundamental group. The quotient group reflects the topological properties of the arithmetic scheme.

Particularly, in \cite{An4} we will use quasi-galois closed covers to give the computation of the \'{e}tale fundamental group of an arithmetic scheme.

\subsection{Outline of the Proof for the Main Theorem}

The whole of \S 3 will be devoted to the proof of the Main Theorem
of the present paper, where we will proceed in several subsections.

In \S 3.1 we will recall some preliminary facts on affine structures
on arithmetic schemes (see [1]). Here, affine structures on a scheme
behave like differential structures on a differential manifold. In
\S 3.2 we will prove that a quasi-galois closed arithmetic variety
has one and only one maximal affine structure among others with
values in a fixed algebraically closed field (see \emph{Proposition
3.9}).

In \S 3.3 we will define conjugations of a given field and a
quasi-galois extension of a field in an evident manner. For the case
of algebraic extensions, \textquotedblleft
conjugation\textquotedblright\ is exactly \textquotedblleft
conjugate\textquotedblright\ and \textquotedblleft
quasi-galois\textquotedblright\ is exactly \textquotedblleft
normal\textquotedblright.

Let $K$ be a finitely generated extension
 of a fixed field $k$.  In \S 3.4 we will demonstrate that $K$
is quasi-galois over $k$ if and only if $K$ has only one conjugation
over $k$ (see \emph{Corollary 3.14}). Moreover, $K$ is a Galois
extension of $k$ if $K$ is quasi-galois and separably generated over
$k$ (see the proof of \emph{Theorem 3.26}).

Then conjugations and quasi-galois extensions for fields will be
geometrically realized in arithmetic varieties. In \S 3.5 we will
define conjugations of an open subset in an arithmetic variety in an
evident manner. An open subset of an arithmetic variety is said to
have a quasi-galois set of conjugations if all of its conjugations
can be affinely realized in the variety.

Now let $\phi :X \rightarrow Y$ be a surjective morphism of finite
type between arithmetic varieties. Suppose that $X$ is quasi-galois
closed over $Y$ by the structure morphism $\phi$.

In \S 3.6 we will establish a relationship between the conjugations
of fields and the conjugations of open subsets in arithmetic
varieties. In deed, the discussions on fields and schemes are
parallel. It will be proved that affine open sets in $X$ have
quasi-galois sets of conjugations (see \emph{Theorem 3.23}) and that
 the function field $k(X)$ is canonically a quasi-galois
extension of the function field $k(Y)$ (see \emph{Theorem 3.24}).

In \S 3.7 we will prove that the automorphism group of $X$ over $Y$
is isomorphic the Galois group of the function field $k(X)$ over
$k(Y)$ (see \emph{Theorem 3.26}), which is the dominant part of the
Main Theorem in the paper. Finally in \S 3.8 we will complete the
proof for the Main Theorem of the paper.

\subsection*{Acknowledgements}

The author would like to express his sincere gratitude to Professor
Li Banghe for his advice and instructions on algebraic
geometry and topology. Thanks for an anonymous referee's comments.

\section{Notation and Definitions}

\subsection{Convention}

In this paper, an \textbf{arithmetic variety} is an integral scheme
of finite type over $ Spec\left( \mathbb{Z}\right) $. A
$k-$\textbf{variety} is an
 integral scheme of finite type over a field $k$.
By a \textbf{variety}, we will understand an arithmetic variety or a
$k-$variety. Let $k(X)\triangleq \mathcal{O}_{X,\xi}$ denote the
function field of a variety $X$ (with generic point $\xi$).

Let $X$ and $Y$ be varieties over a fixed variety $Z$.  $Y$ is
said to be a \textbf{conjugate} of $X$ over $Z$ if there is an
isomorphism $\sigma :X\rightarrow Y$ over $Z$. Let
$Aut\left( X/Z\right)$ denote the group of automorphisms of $X$ over
$Z$.

Let $D$ be an integral domain. Denote by $Fr(D)$ the field of
fractions on $ D $. If $D$ is a subring of a field
$\Omega$, $Fr(D)$ will be assumed to be contained in $\Omega$.

Let $E$ be an extension of a field $F$. Note that here
 $E$ is not necessarily algebraic over $F$.  Recall that $E$ is a \textbf{Galois extension} of $F$ if
$F$ is the invariant subfield of the Galois group $Gal(E/F)$.

\subsection{Affine covering with values in a given field}

Let $(X,\mathcal{O} _{X})$ be a scheme. As usual, an
affine covering of the scheme $(X,\mathcal{O} _{X})$ is a
family $\mathcal{C}_{X}=\{(U_{\alpha },\phi _{\alpha };A_{\alpha
})\}_{\alpha \in \Delta }$ such that for each $\alpha \in \Delta $,
 $\phi _{\alpha }$ is an isomorphism from an open set $U_{\alpha
}$ of $X$ onto the spectrum $Spec{A_{\alpha }}$ of a
commutative ring $A_{\alpha }$. Each $(U_{\alpha
},\phi _{\alpha };A_{\alpha })\in \mathcal{C}_{X}$ is called a
\textbf{local chart}. An affine covering
$\mathcal{C}_{X}$ of
$(X, \mathcal{O}_{X})$ is said to be \textbf{reduced} if
$U_{\alpha}\neq U_{\beta} $ holds for any $\alpha\neq \beta$ in
$\Delta$.

Sometimes, we will denote by $(X,\mathcal{O}_{X};\mathcal{C}_{X})$ a scheme
$(X,\mathcal{O}_{X})$  with a given affine covering
$\mathcal{C}_{X}$. For the sake of brevity, a local chart
$(U_{\alpha},\phi_{\alpha};A_{\alpha })$  will be denoted
by $U_{\alpha}$ or $(U_{\alpha},\phi_{\alpha})$.

Let $\mathfrak{Comm}$ be the category of commutative rings with
identity. Fixed a subcategory $\mathfrak{Comm}_{0}$ of
$\mathfrak{Comm}$. An affine covering
$\{(U_{\alpha},\phi_{\alpha};A_{\alpha })\}_{\alpha \in \Delta}$ of
$(X, \mathcal{O}_{X})$ is said to be \textbf{with values} in
$\mathfrak{Comm}_{0}$ if $\mathcal{O}_{X}(U_{\alpha})=A_{\alpha}$
holds and $A_{\alpha }$ is contained in $\mathfrak{Comm}_{0}$ for
each $\alpha \in \Delta$.

In particular, let $\Omega$ be a field (large enough) and let
$\mathfrak{Comm}(\Omega)$ be the category consisting of the subrings
of $\Omega$ and their isomorphisms. An affine covering
$\mathcal{C}_{X}$ of $(X, \mathcal{O}_{X})$ with values in
$\mathfrak{Comm}(\Omega)$ is said to be \textbf{with values in the
field $\Omega$}.

\subsection{Quasi-galois closed varieties}

Let $X$ and $Y$ be two varieties and let $f:X\rightarrow Y$ be a
surjective morphism of finite type.

\begin{definition}
The variety $X$ is said to be \textbf{quasi-galois closed}  over $Y$
by $f$ if there is an algebraic closed field $\Omega$ and a reduced
affine covering $\mathcal{C}_{X}$ of $X$ with values in $\Omega$
such that for any conjugate $Z$ of $X$ over $Y$  the following
conditions are satisfied:

$(i)$ $(X,\mathcal{O}_{X})=(Z,\mathcal{O}_{Z})$ holds if
$(Z,\mathcal{O}_{Z})$ has a reduced affine coverings with values in
$\Omega$.

$(ii)$ Each local chart contained in $\mathcal{C}_{Z}$ is contained
in $ \mathcal{C}_{X}$ for any reduced affine covering
$\mathcal{C}_{Z}$ of $(Z,\mathcal{O}_{Z})$ with values in $\Omega$.
\end{definition}

In particular, if $Y$ is $Spec(\mathbb{Z})$ or $Spec(k)$, such a
variety $X$ is said to be a \textbf{quasi-galois closed variety}.

\begin{remark}
\textbf{The existence of quasi-galois closed varieties.}

$(i)$ For the case of varieties, the finite group actions on
varieties can produce quasi-galois closed varieties (For example, see \cite{EGA,SGA1,GIT,VS1,SV2}).

$(ii)$ For the case of schemes, there is another way to obtain
quasi-galois closed schemes. Let $X$ be a scheme with a finite
number of conjugates. Then the disjoint union of the conjugates of
$X$ will be quasi-galois closed over $X$.

$(iii)$ For a general case, in \cite{An3} we will prove the existence of quasi-galois closed schemes over arithmetic schemes.
\end{remark}

\section{Statement of The Main Theorem}

Here is the \emph{Main Theorem} of the present paper, which will be proved
in \S 3.

\begin{theorem}
\textbf{\emph{(Main Theorem).}} Let $X$ and $Y$ be two arithmetic
varieties. Assume that $X$ is quasi-galois closed over $Y$ by a
surjective morphism $\phi$ of finite type. Then there are the following statements.

\begin{itemize}
\item $f$ is affine.

\item $k\left( X\right) $ is canonically a Galois
extension of $k(Y)$.

\item There is a group isomorphism
$$
{Aut}\left( X/Y\right) \cong Gal(k\left( X\right) /k(Y)).
$$

\item Particularly, let $\dim X=\dim Y$.
Then  $X$ is a pseudo-galois cover of  $Y$ in the sense of Suslin-Voevodsky.
\end{itemize}
\end{theorem}

\begin{remark}
By the first property in \emph{Theorem 2.1} it is seen that
there exists a
nice relationship between quasi-galois closed arithmetic varieties
and Galois extensions of functions fields in several variables.
\end{remark}

\begin{remark}
If $\dim X=\dim Y$, it is seen that quasi-galois closed arithmetic varieties behave like
Galois extensions of number fields and their automorphism groups can
be regarded as the Galois groups of the field extensions. If $\dim X>\dim Y$,
\emph{Theorem 2.1} can be regarded as a
generalization of that in \emph{Proposition 1.1} in \cite{SGA1}, \emph{Page
106} for function fields in several variables.
\end{remark}

\begin{remark}
We have attempted to use the data of such varieties $X/Y$ to describe a
given finite Galois extension $E/F$ in such a manner that $X/Y$ are said to be a \emph{model} for
$E/F$ if the Galois group $Gal\left( E/F\right) $ is isomorphic to
the automorphism group ${Aut}\left( X/Y\right)$ (for example, see \cite{SGA1,GIT,Raskind,VS1,SV2}).
Hence, \emph{Theorem 2.1} afford us such a model for function fields in several variables.
\end{remark}

\begin{remark}
If $\dim X=\dim Y$, we have pseudo-galois covers of arithmetic varieties in the sense of Suslin-Voevodsky (see \cite{VS1,SV2}); it is seen that there is no essential difference between our \textquotedblleft quasi-galois closed\textquotedblright and \textquotedblleft pseudo-galois cover\textquotedblright.
However, suppose $\dim X>\dim Y$. Then it is seen that there is a main difference between the two types of covers. For example,  $Spec(\mathbb{Z}[t])/Spec(\mathbb{Z})$ is quasi-galois closed but not pseudo-galois, where $t$ is a variable over $\mathbb{Q}$.
Hence, \emph{Theorem 2.1} gives us a sufficient
condition for the existence of such a pseudo-galois cover in a more
generalized case in the category of arithmetic varieties, where the function fields are in several variables.
\end{remark}

\section{Proof of the Main Theorem}

In this section we will proceed in several subsections to prove
 the main theorem of the paper.

\subsection{Affine structures}

Let us recall some preliminary results on affine structures (see
\cite{An}) which will be used in the following subsections. Here, affine
structures on a schemes can be regarded as a counterpart of
differential structures on a manifold in topology (for example, see \cite{STR}).

Let $\mathfrak{Comm}$ be the category of commutative rings with
identity. Fixed a subcategory $\mathfrak{Comm}_{0}$ of
$\mathfrak{Comm}$.

\begin{definition}
A \textbf{pseudogroup $\Gamma$ of affine transformations}
(\textbf{with values in $\mathfrak{Comm}_{0}$}) is a subcategory of
$\mathfrak{Comm}_{0}$ such that the algebra isomorphisms contained
in $\Gamma$ satisfying the conditions $(i)-(v)$:

$\left( i\right) $ Each $\sigma\in\Gamma$ is an isomorphism between
algebras $dom\left( \sigma\right) $ and $rang\left( \sigma\right)$
contained in $\mathfrak{Comm}_{0}$, called the \textbf{domain} and
\textbf{range} of $\sigma$, respectively.

$\left( ii\right) $ Let $\sigma\in\Gamma$. Then the inverse
$\sigma^{-1}$ is contained in $\Gamma.$

$\left( iii\right) $ The identity map $id_{A}$ on $A$ is contained in $%
\Gamma$ if there is some $%
\delta\in\Gamma$ with $dom\left( \delta\right) =A.$

$\left( iv\right) $ Let $\sigma\in\Gamma$. Then the isomorphism
induced by $\sigma$ defined on the localization $dom\left(
\sigma\right)_{f}$ of the algebra $dom\left( \sigma\right)$ at any
nonzero $f\in dom\left( \sigma\right) $ is contained in $\Gamma.$

$\left( v\right) $ Let $\sigma,\delta\in\Gamma$. Assume for some
$\tau\in\Gamma$ there are isomorphisms $dom\left( \tau\right) \cong
dom\left( \sigma\right) _{f}$ and $dom\left( \tau\right) \cong
rang\left( \delta\right) _{g}$ with $0\not= f\in dom\left(
\sigma\right) $ and $0\not= g\in rang\left( \delta\right) .$ Then
the isomorphism factorized by $dom\left( \tau\right) $ from
$dom\left( \sigma\right) _{f}$ onto $rang\left( \delta\right) _{g}$
is contained in $\Gamma$.
\end{definition}

Let $X$ be a topological space and let $\Gamma$ be a pseudogroup of
affine transformations with values in $\mathfrak{Comm}_{0}$.

\begin{definition}
An \textbf{affine} $\Gamma-$\textbf{atlas $\mathcal{A}$ on $X$}
(\textbf{with values in $\mathfrak{Comm}_{0}$}) is a collection of
triples $\left( U_{j},\varphi _{j}; A_{j}\right) $ with $j\in
\Delta$, called \textbf{local charts}, satisfying the conditions
$(i)-(iii)$:

$\left( i\right) $ For every $\left( U_{j},\phi_{j}; A_{j}\right)
\in \mathcal{A}$, $U_{j}$ is an open subset of $X$ and $\phi _{j}$
is an homeomorphism of $U_{j}$ onto $Spec\left( A_{j}\right)$ with
$A_{j}\in\Gamma$ such that $U_{i}\not= U_{j}$ holds for any $i\not=
j$ in $\Delta$.

For the sake of brevity, such a triple $\left( U_{j},\phi_{j};
A_{j}\right)$ will be denoted sometimes by $U_{j}$ or by a pair
$\left( U_{j},\phi_{j}\right)$.

{$\left( ii\right) $ $\bigcup_{j\in\Delta} U_{j}$ is an open
covering of $X.$}

{$\left( iii\right) $ Take any $\left( U_{i},\phi_{i},A_{i}\right)
,\left( U_{j},\phi_{j},A_{j}\right) \in \mathcal{A}$ with $U_{i}\cap
U_{j}\not =\varnothing$. Then there is a local chart $\left(
W_{ij},\phi _{ij}\right) \in \mathcal{A}$ with $ W_{ij}\subseteq
U_{i}\cap U_{j}$ such that the isomorphism between the localizations
$\left( A_{j}\right) _{f_{j}}$ and $\left( A_{i}\right) _{f_{i}}$
induced by the map $
\phi_{j}\circ\phi_{i}^{-1}\mid_{W_{ij}}:\phi_{i}(W_{ij})\rightarrow%
\phi_{j}(W_{ij}) $ is contained in $\Gamma$, where $\phi_{i}\left(
W_{ij}\right) \cong Spec\left( A_{i}\right) _{f_{i}}$ and
$\phi_{j}\left( W_{ij}\right) \cong Spec\left( A_{j}\right)
_{f_{j}}$ are homeomorphic
 for some $f_{i}\in A_{i}$ \text{ and }$f_{j}\in A_{j}.
$}
\end{definition}

\begin{definition}
{Two affine $\Gamma-$atlases $\mathcal{A}$ and $\mathcal{A}%
^{\prime} $ on $X$\ are said to be $\Gamma-$\textbf{compatible} if
the condition below is satisfied:}

{Take any $\left( U,\phi, A\right) \in \mathcal{A}$ and $\left(
U^{\prime },\phi^{\prime},A^{\prime}\right) \in
\mathcal{A}^{\prime}$ with $U\cap U^{\prime}\not =\varnothing.$ Then
there is a local chart $\left( W,\phi ^{\prime\prime}\right) \in
\mathcal{A}\bigcap \mathcal{A}^{\prime}$ with $W\subseteq U\cap
U^{\prime}$ such that the isomorphism between the localizations $A
_{f}$ and $\left( A^{\prime}\right)
_{f^{\prime}}$ induced by the map $\phi^{\prime}\circ\phi^{-1}%
\mid_{W}:\phi(W)\rightarrow \phi^{\prime}(W)$ is contained in
$\Gamma$, where $\phi\left( W\right) \cong SpecA _{f}$ and
$\phi^{\prime}\left( W\right) \cong Spec\left( A^{\prime}\right)
_{f^{\prime}}$ are homeomorphic for some $f\in A$ and $f^{\prime}\in
A^{\prime}.$}
\end{definition}

By an \textbf{affine }$\Gamma-$\textbf{structure} on $X$ (with
values in $\mathfrak{Comm}_{0}$)
 we understand a maximal affine $\Gamma-$atlas $\mathcal{A}\left( \Gamma\right) $
 on $X$. Here, an affine
$\Gamma-$atlas $\mathcal{A}$ on $X$ is said to be \textbf{maximal}
(or \textbf{complete}) if it can not be contained properly in any
other affine $\Gamma-$atlas of $X.$

\begin{remark}
 Fixed a pseudogroup
$\Gamma$ of affine transformations. By Zorn's Lemma it is seen that
for any given affine $\Gamma-$atlas $\mathcal{A}$ on $X$ there is a
unique affine $\Gamma-$structure $\mathcal{A}_{m}$ on $X$ satisfying

$\left( i\right) $ $\mathcal{A}\subseteq \mathcal{A}_{m};$

$\left( ii\right) $ $\mathcal{A}$ and $\mathcal{A}_{m}$ are $\Gamma-$%
compatible.

In such a case, $\mathcal{A}$ is said to be a \textbf{base} for $\mathcal{A}_{m}$ and $%
\mathcal{A}_{m}$ is the affine $\Gamma-$structure defined by $%
\mathcal{A}.$
\end{remark}

\begin{definition}
Let $\mathcal{A}\left( \Gamma\right) $ be a affine
$\Gamma-$structure on $X$. Assume that there is a sheaf
$\mathcal{F}$ of rings on $X$ such that $\left( X,\mathcal{F}\right)
$ is a locally ringed space and that $
\phi_{\alpha\ast}\mathcal{F}\mid_{U_{\alpha}}\left(
SpecA_{\alpha}\right) =A_{\alpha}$ holds for each $\left( U_{\alpha
},\phi_{\alpha};A_{\alpha}\right) \in \mathcal{A}\left(
\Gamma\right) $.

Then $\mathcal{A}\left( \Gamma\right) $ is said to be
\textbf{admissible} on $X$ and $\mathcal{F}$ is said to be an
\textbf{extension}
 of $
\mathcal{A}\left( \Gamma\right) $.
\end{definition}

It is evident that such a sheaf $\mathcal{F}$ on $X$ affords us a
scheme $(X,\mathcal{F})$. That is, an extension of an affine
structure on a space is a scheme.

Let $(X,\mathcal{O}_{X})$ be a scheme and $U_{\alpha}$ an affine
open set of $X$. Take an isomorphism
$(\phi_{\alpha},{\phi_{\alpha}}^{\sharp}):(U_{\alpha},{\mathcal{O}_{X}}_
{\mid U_{\alpha}})\rightarrow
(SpecA_{\alpha},\mathcal{O}_{SpecA_{\alpha}})$. In general, the ring
${\mathcal{O}_{X}}(U)$ is isomorphic to $A_{\alpha}$ by
${\phi_{\alpha}}^{\sharp}$. Here, we choose the ring $A_{\alpha}$ to
be such that ${\mathcal{F}}(U)=A_{\alpha}$ in the definition above.
This can be done according to the preliminary facts on affine
schemes (see \cite{EGA}).

\begin{remark}
It is easily seen that all extensions of a fixed admissible affine
structure on a space are isomorphic schemes (see \cite{An}).
\end{remark}

\subsection{A quasi-galois closed variety has only one maximal
affine structure among others
 with values in a fixed field}

Let $\mathfrak{Comm}_{/k}$ be the category of finitely generated
algebras (with identities) over a given field $k$. We will consider
the pseudogroup of affine transformations with values in
$\mathfrak{Comm}_{/k}$ in this subsection.

Fixed a $k-$variety $\left( X,\mathcal{O}_{X}\right) $ with a given
reduced affine covering $ \mathcal{C}_{X}$. That is, each reduced
affine covering gives us a pseudogroup of affine transformations in
a natural manner.

In deed, define $\Gamma (\mathcal{C}_{X})$ to be the set
of identities $1_{A_{\alpha}}:A_{%
\alpha}\rightarrow A_{\alpha}$  and isomorphisms
$\sigma_{\alpha\beta}:\left( A_{\alpha}\right)
_{f_{\alpha}}\rightarrow \left( A_{\beta}\right) _{f_{\beta}} $ of
$k-$algebras for any $0\not= f_{\alpha }\in A_{\alpha }$ and $0\not=
f_{\beta }\in A_{\beta }$, where $A_{\alpha }$ and $A_{\beta }$ are
contained in $\mathfrak{Comm}/k$ such that there are some affine
open subsets $U_{\alpha }$ and $U_{\beta }$ of $X$ with $(U_{\alpha
},\phi _{\alpha }),(U_{\beta }, \phi _{\beta })\in \mathcal{C}_{X}$
satisfying $$\phi _{\alpha }\left( U_{\alpha }\right) =SpecA_{\alpha
},\phi _{\beta }\left( U_{\beta }\right) =SpecA_{\beta }.$$

Then $\Gamma (\mathcal{C}_{X})$ is a pseudogroup in $\mathfrak{%
Comm}_{/k}$, called the (\textbf{maximal}) \textbf{pseudogroup of
affine transformations} in $\left(
X,\mathcal{O}_{X};\mathcal{C}_{X}\right) $.

It is seen that $\mathcal{C}_{X})$ is an affine $\Gamma
(\mathcal{C}_{X})-$atlas on $X$. Denote by $\mathcal{A}
(\mathcal{C}_{X})$ the affine $\Gamma (\mathcal{C}_{X})-$structure
on $X$ defined by an affine $\Gamma (\mathcal{C}_{X})-$atlas
$\mathcal{C}_{X}$ on $X$.

\begin{remark}
Let $\mathfrak{Comm}_{/\mathbb{Z}}$ be the category of finitely
generated algebras (with identities) over $\mathbb{Z}$. We can
similarly define a pseudogroup $\Gamma$ of affine transformations
with values in $\mathfrak{Comm}_{/\mathbb{Z}}$ and then discuss
 affine $\Gamma-$structures such as the above.

 In particular, let $\Omega$ be a field and let
$\mathfrak{Comm}(\Omega)$ be the category consisting of the subrings
of $\Omega$ and their isomorphisms. An affine structure
 of a variety $X$ with values in
$\mathfrak{Comm}(\Omega)$ is said to be \textbf{with values in the
field $\Omega$}.
\end{remark}

\begin{remark} Let $\left( X,\mathcal{O}_{X}\right) $ be a variety
(i.e, an arithmetic variety or a $k-$variety).

$(i)$ Different reduced affine coverings $ \mathcal{C}_{X}$ have
different pseudogroups $\Gamma ( \mathcal{C}_{X})$ of affine
transformations.

$(ii)$ Each reduced affine covering $\mathcal{C}_{X}$ is an
admissible affine atlas. In particular, $\mathcal{O}_{X}$ is an
extension of $\mathcal{C}_{X}$ on the underlying space $X$.

$(iii)$ As an admissible affine atlas, each reduced affine covering
$\mathcal{C}_{X}$ can have many extensions $\mathcal{F}_{X}$ on the
space $X$. With such an extension we have a scheme
$(X,\mathcal{F}_{X})$, called an \textbf{associate scheme} of
$\left( X,\mathcal{O}_{X}\right) $.

$(iv)$ By \emph{Remark 3.6} it is seen that all associate schemes of
$\left( X,\mathcal{O}_{X}\right) $ are isomorphic.
\end{remark}

\begin{proposition}
Let $X$ and $Y$ be varieties with $X$ quasi-galois closed over $Y$.
Then the variety $X$ has one and only one affine structure
$\mathcal{A}(\mathcal{O}_{X})$ with values in an algebraic closed
field $\Omega$ satisfying the below properties:

$(i)$ The structure sheaf $\mathcal{O}_{X}$ is an extension of
$\mathcal{A}(\mathcal{O}_{X})$.

$(ii)$ By set inclusion, $\mathcal{A}(\mathcal{O}_{X})$ is maximal
among the whole of the affine structures on the underlying space of
$X$ with values in $\Omega$. That is, take any affine
$\Gamma-$structure $\mathcal{B}$ on the space of $X$ with values in
$\Omega$. Then we must have $\mathcal{B}\subseteq
\mathcal{A}(\mathcal{O}_{X})$ and $\Gamma \subseteq
\Gamma(\mathcal{A}(\mathcal{O}_{X}))$, where
$\Gamma(\mathcal{A}(\mathcal{O}_{X}))$ is the maximal pseudogroup of
affine transformations of
$(X,\mathcal{O}_{X};\mathcal{A}(\mathcal{O}_{X}))$.

In particular, we can choose $\Omega$ to be an fixed algebraic
closure of the function field $k(X)$ of $X$.

Here, $\mathcal{A}(\mathcal{O}_{X})$ will be called the
\textbf{natural affine structure} of $(X,\mathcal{O}_{X})$ with
values in $\Omega$.
\end{proposition}

\begin{proof}
Let $\mathcal{C}_{X}$ and $\mathcal{C}_{X}^{\prime}$ be two affine
structures on the underlying space $X$ with respect to pseudogroups
$\Gamma$ and $\Gamma^{\prime}$ respectively, which are both with
values in some field $\Omega$. As an affine structure is a maximal
affine atlas, it is clear that $\mathcal{C}_{X}$ and
$\mathcal{C}_{X}^{\prime}$ are reduced affine coverings on the space
$X$.  By \emph{Definition 1.1}, we must have either
$$\mathcal{C}_{X}\subseteq \mathcal{C}_{X}^{\prime},\Gamma\subseteq
\Gamma^{\prime}$$ or $$\mathcal{C}_{X}\supseteq
\mathcal{C}_{X}^{\prime},\Gamma\supseteq \Gamma^{\prime}.$$

Let $\Sigma$ be the set of affine structures on the underlying space
$X$ with values in $\Omega$. By set inclusion, $\Sigma$ is a
partially ordered set since any two affine structures are compatible
with the pseudogroups of affine transformations.

Hence, $\Sigma$ is totally ordered. The unique maximal element in
$\Sigma$ is the desired affine structure, where we choose the field
$\Omega$ to be an algebraic closure of the functional field
$\mathcal{O}_{X,\xi}=k(X)$ of  $X$.
\end{proof}

\subsection{Definition for conjugations of a given field}

Let $K$ be an extension of a field $k$. Here $K/k$ is not
necessarily algebraic. $K$ is said to be \textbf{$k-$quasi-galois}
(or, \textbf{quasi-galois} over $k$) if each irreducible polynomial
$f(X)\in F[X]$ that has a root in $K$ factors completely in $K\left[
X\right] $ into linear factors for any intermediate field
$k\subseteq F\subseteq K$.

Let $E$ be a finitely generated extension of $k$. The elements
\begin{equation*}
w_{1},w_{2},\cdots ,w_{n}\in E\setminus k
\end{equation*}
are said to be a \textbf{$(r,n)-$nice $k-$basis} of $E$ (or simply, a
\textbf{nice $k-$basis}) if the following conditions are satisfied:

$E=k(w_{1},w_{2},\cdots ,w_{n})$;

$w_{1},w_{2},\cdots ,w_{r}$ constitute a transcendental basis of $E$ over $k$%
;

$w_{r+1},w_{r+2},\cdots ,w_{n}$ are linearly independent over $%
k(w_{1},w_{2},\cdots ,w_{r})$, where $0\leq r\leq n$.

\begin{definition}
Let $E$ and $F$ be two finitely generated extensions of a field $k$.
$F$ is said to be a \textbf{$k-$conjugation} of $E$ (or, a
\textbf{conjugation} of $E$ over $k$) if there is a $(r,n)-$nice
$k-$basis $w_{1},w_{2},\cdots ,w_{n} $ of $E$ such that $F$ is a
conjugate of $E$ over $k(w_{1},w_{2},\cdots ,w_{r})$.

We will denote by $\tau_{(r,n)}$ such an isomorphism from $F$ onto
$E$ over $k(w_{1},w_{2},\cdots ,w_{r})$ with respect to the
$(r,n)-$nice $k-$basis.
\end{definition}

\begin{remark}
Let $F$ be a $k-$conjugation of $E$. Then $F$ is contained in the
algebraic closure $\overline{E}$ of $E$.

It will be proved that a finitely generated field is quasi-galois if
and only if it has only one conjugation (see \emph{Corollary 3.14}).
For the case of algebraic extensions, this is exactly to say that a
field is normal if and only if it has only one conjugate field.
\end{remark}

\subsection{A quasi-galois field has only one conjugation}

We give the below criterion for  a quasi-galois field by
conjugations, which behaves like a normal field and its conjugate
field for the case of algebraic extensions.

\begin{theorem}
Let $K$ be a finitely generated extension of a field $k$. The
following statements are equivalent.

$\left( i\right) $ $K$ is a quasi-galois field over $k$.

$\left( ii\right) $ Take any $x\in K$ and any subfield $k\subseteq
F\subseteq K$. Then each conjugation of $ F\left( x\right) $ over
$F$ is contained in $K$.

$\left( iii\right) $ Each $k-$conjugation of $K$ is contained in
$K$.
\end{theorem}

\begin{proof}
$\left( i\right) \implies \left( ii\right) .$ Take any $x\in K$ and
any $k\subseteq F\subseteq K.$ If $x$ is a variable over $F$, the
field $F\left( x\right) $ is the unique $k-$conjugation of $F\left(
x\right) $ in $\overline{F\left( x\right) }$
($\subseteq\overline{K}$). If $x$ is algebraic over $F$, a
$F-$conjugation of $F\left( x\right)$ which is exactly a
$F-$conjugate of $F\left( x\right)$ is contained in $K$ by the
assumption that $K$ is $k-$quasi-galois; then all $ F-$conjugates of
$F\left( x\right)$ in $\overline{F\left( x\right) }$
($\subseteq\overline{K}$) is contained in $K$.

$\left( ii\right) \implies \left( i\right) .$ Let $F$ be a field
with $ k\subseteq F\subseteq K$. Take any irreducible polynomial
$f\left( X\right) $ over $F.$ Suppose that $x\in K$ satisfies the
equation $f\left( x\right) =0$. Then such an $F-$conjugation of
$F(x)$ is an $F-$conjugate. By $\left( ii\right) $ it is seen that
every $F-$conjugate $z\in \overline{F}$ of $x $ is contained in $K$;
hence, $K$ is quasi-galois over $k.$

$\left( ii\right) \implies \left( iii\right) .$ Hypothesize that
there is a $ k-$conjugation $H$ of $K$ in $\overline{K}$ is not
contained in $K,$ that is, $H\setminus K$ is a nonempty set. Take
any $x_{0}\in H\setminus K$.

Choose a $(r,n)-$nice $k-$basis $w_{1},w_{2},\cdots ,w_{n} $ of $K$
which make $H$ be a $k-$conjugation of $K$. By \emph{Remark 3.11} it
is seen that $H$ is contained in the algebraic closure of
$k(w_{1},w_{2},\cdots ,w_{n} )$. As $w_{1},w_{2},\cdots ,w_{r} $ are
all variables over $k$, it is seen that $w_{1},w_{2},\cdots ,w_{r} $
are all contained in the intersection of $H$ and $K$. By
\emph{Definition 3.10} it is seen that there is an isomorphism
$\sigma:H\rightarrow K$ of fields over $k(w_{1},w_{2},\cdots ,w_{r}
)$.

It is evident that the specified element $x_{0}$ must be algebraic
over $k(w_{1},w_{2},\cdots ,w_{r} )$. Then the field
$k(w_{1},w_{2},\cdots ,w_{r} ,x_{0})$ is a conjugate of the field
$k(w_{1},w_{2},\cdots ,w_{r} ,\sigma(x_{0}))$ over
 $k(w_{1},w_{2},\cdots ,w_{r} )$.

From $\left( ii\right) $ we have $k(w_{1},w_{2},\cdots ,w_{r}
,x_{0})\subseteq K$. In particular, $x_{0}$ is contained in $K$,
which is in contradiction with the hypothesis above. Therefore,
every $k-$conjugation of $K$ is in $K.$

$\left( iii\right) \implies \left( ii\right) .$ Take any $x\in K$
and any field $F$ such that $k\subseteq F\subseteq K$. If $x$ is a
variable over $F,$
$F\left( x\right) $ is the unique $F-$conjugation in $\overline{K}$ of $%
F\left( x\right) $ itself by \emph{Remark 3.11} again; hence,
$F\left( x\right) $ is contained in $K.$

Now suppose that $x$ is algebraic over $F.$ Let $z\in \overline{K}$
be an $F-$conjugate of $x$. If $F=K,$ we have $\sigma_{x} =id_{K};$
then $z=x\in K$. If $F\not=K$, from \emph{Lemma 3.13} below we have
a field $F\left( z,v_{1},v_{2},\cdots ,v_{s},w_{s+1},\cdots,
w_{m}\right)$ that is an $F-$conjugation of $K$;
 it is seen that the
element $z$ is contained in an $F-$conjugation of $K$; as $k
\subseteq K$, an $F-$conjugation of $K$ must be an $k-$conjugation
of $K$; by $(iii)$ we must have $z\in K$. This proves $(ii)$.
\end{proof}

\begin{lemma}
Fixed a finitely generated extension $K$ of a field $k$ and a field
$F$ with $k\subseteq F\subsetneqq K$. Let $x\in K$ be algebraic over
$F$ and let $z$ be a conjugate of $x$ over $F$. Then there is a
$(s,m)-$nice $F\left( x\right) -$basis $v_{1},v_{2},\cdots ,v_{m}$
of $K$ and an $F-$isomorphism $\tau$ from the field
$$K=F\left( x,v_{1},v_{2},\cdots ,v_{s},v_{s+1},\cdots, v_{m}\right)$$
onto a field of the form $$F\left( z,v_{1},v_{2},\cdots
,v_{s},w_{s+1},\cdots, w_{m}\right)$$ such that $$\tau (x)=z,\tau
(v_{1})=v_{1},\cdots,\tau (v_{s})=v_{s}$$ where
$w_{s+1},w_{s+2},\cdots, w_{m}$ are elements contained in an
extension of $F$. In particular, we have
$$w_{s+1}=v_{s+1},w_{s+2}=v_{s+2},\cdots, w_{m}=v_{m}$$ if $z$ is not contained in
$F(v_{1},v_{2},\cdots,v_{m})$.
\end{lemma}

\begin{proof}
We will proceed in two steps according to the assumption that $s=0$
or $s\not= 0$.

\emph{Step 1}. Let $s\not= 0$.
That is, $v_{1}$ is a variable over $F\left( x\right) $.

Let $\sigma_{x}$ be the $F-$isomorphism between fields $F(x)$ and
$F(z)$ with $\sigma_{x}(x)=z$. From the isomorphism $\sigma_{x}$ we
obtain an isomorphism $\sigma _{1}$ of $F\left( x,v_{1}\right) $
onto $F\left( z,v_{1}\right) $ defined by
$$
\sigma_{1}:\frac{f(v_{1})}{g(v_{1})}\mapsto \frac{\sigma _{x}(f)(v_{1})}{%
\sigma _{x}(g)(v_{1})}
$$
for any polynomials $ f[X_{1}],g[X_{1}]\in F\left( x\right) [X_{1}]
$ with $g[X_{1}]\neq 0$.

It is easily seen that $g(v_{1})=0$ if and only if ${\sigma _{x}(g)(v_{1})}%
=0 $. Hence, the map $\sigma_{1}$ is well-defined.

Similarly, for the elements $v_{1},v_{2},\cdots ,v_{s}\in K$ that
are variables over $F(x)$, there is an isomorphism
$$
\sigma _{s}: F\left( x,v_{1},v_{2},\cdots ,v_{s}\right)\longrightarrow
F\left( z,v_{1},v_{2},\cdots ,v_{s}\right)
$$
of fields defined by
$$
x\longmapsto z\text{ and }v_{i}\longmapsto v_{i}
$$
for $1\leq i\leq s$, where we have the restrictions
$$
\sigma_{i+1}|_{F\left( x,v_{1},v_{2},\cdots ,v_{i}\right)}=\sigma_{i}
.$$

If $s=m$, we have $K=F\left( x,v_{1},v_{2},\cdots ,v_{s}\right)$ and
it follows that the field $F\left( z,v_{1},v_{2},\cdots
,v_{s}\right)$ is an $F-$conjugation of $K$. We put $s\leqslant
m-1$.

\emph{Step 2}. Let $s=0$. That is exactly to consider the case
$v_{s+1}\in K$ since $v_{s+1}$ is algebraic over the field
$F(v_{1},v_{2},\cdots ,v_{s})\subseteq K$. We have two cases for the
element $v_{s+1}$.

\emph{Case (i)}. Suppose that $z$ is not contained in
$F(v_{1},v_{2},\cdots,v_{s+1})$.

We have an isomorphism $\sigma _{s+1}$ between the fields $F\left(
x,v_{1},v_{2},\cdots ,v_{s+1}\right) $ and $F\left(
z,v_{1},v_{2},\cdots ,v_{s+1}\right) $ given by
\begin{equation*}
x\longmapsto z\text{ and }v_{i}\longmapsto v_{i}
\end{equation*}
with $1\leq i\leq s+1.$

The map $\sigma_{s+1}$ is well-defined. In deed, by the below
\textbf{Claim$^{\dag}$} it is seen that $f\left( v_{s+1}\right) =0$
holds if and only if $\sigma_{s}\left( f\right) \left(
v_{s+1}\right) =0$ holds for any polynomial $f\left( X_{s+1}\right)
\in F\left( x,v_{1},v_{2},\cdots,v_{s}\right) \left[ X_{s+1}\right]
$.

\emph{Case (ii)}. Suppose that $z$ is contained in the field
$F(v_{1},v_{2},\cdots,v_{s+1})$.

By the below \textbf{Claim$^{\dag\dag}$} we have an  element
$v_{s+1}^{\prime}$ contained in an extension of $F$ such that the
fields $F(x,v_{s+1})$ and $F(z, v_{s+1}^{\prime})$ are isomorphic
over $F$.

Then by the same procedure as in \emph{Case (i)} of
\textbf{Claim$^{\dag}$} it is seen that the fields
$F(x,v_{s+1},v_{1},v_{2},\cdots,v_{s})$ and
$F(z,v_{s+1}^{\prime},v_{1},v_{2},\cdots,v_{s})$ are isomorphic over
$F$.

Hence, in such a manner we have an $F-$isomorphism $ \tau$ from the
field
$$F\left( x,v_{1},v_{2},\cdots ,v_{s},v_{s+1},\cdots, v_{m}\right)$$
onto the field of the form $$F\left( z,v_{1},v_{2},\cdots
,v_{s},w_{s+1},\cdots, w_{m}\right)$$ such that $$\tau (x)=z,\tau
(v_{1})=v_{1},\cdots,\tau (v_{s})=v_{s},$$ where
$w_{s+1},w_{s+2},\cdots, w_{m}$ are elements contained in an
extension of $F$. This completes the proof of the lemma.
\end{proof}

\textbf{Claim$^{\dag}$}. Given any $f\left( X,X_{1},X_{2}
,\cdots,X_{s+1}\right)$ in the polynomial ring $F\left[
X,X_{1},X_{2},\cdots,X_{s+1}\right] $. Suppose that $z$ is not
contained in the field $F(v_{1},v_{2},\cdots,v_{s+1})$. Then
$f\left( x,v_{1},v_{2},\cdots,v_{s+1}\right) =0$ holds if and only
if $f\left( z,v_{1},v_{2},\cdots,v_{s+1}\right) =0$ holds.

\begin{proof}
Here we use Weil's algebraic theory of specializations (See \cite{weil})
 to prove the claim. For $v_{s+1}$ there are two cases:
$$v_{s+1}\in \overline{F};$$
$$v_{s+1}\in \overline{F(v_{1},v_{2},\cdots,v_{s+1})}\setminus
\overline{F},$$ where $\overline{F}$ denotes the algebraic closure
of the field $F$.

\emph{Case (i)}. Let $v_{s+1}\in
\overline{F(v_{1},v_{2},\cdots,v_{s+1})}\setminus \overline{F}$.

By \emph{Theorem 1} in \cite{weil}, \emph{Page 28},  it is clear that
$\left( z\right) $ is a (generic) specialization of $\left( x\right)
$ over $F$ since $z$ and $x$ are conjugates over
 $F$.
From \emph{Proposition 1} in \cite{weil}, \emph{Page 3}, it is seen that
$F\left( v_{1},v_{2},\cdots ,v_{s+1}\right) $ and the field $F(x)$
are free with respect to each other over $F$ since $x$ is algebraic
over $F$. That is, $F\left( v_{1},v_{2},\cdots ,v_{s+1}\right) $ is
a free field over $F$ with respect to $(x)$. By \emph{Proposition 3}
in \cite{weil}, \emph{Page 4}, it is seen that $F\left( v_{1},v_{2},\cdots
,v_{s+1}\right) $ and the algebraic closure $\overline{F}$ are
linearly disjoint over $F$. That is, $F\left( v_{1},v_{2},\cdots
,v_{s+1}\right) $ is a regular extension of $F$ (For detail, see
\cite{weil}, \emph{Page 18}).

Then by \emph{Theorem 5} in \cite{weil}, \emph{Page 29}, it is seen that
$\left( z,v_{1},v_{2},\cdots,v_{s+1}\right) $ is a (generic)
specialization of $\left( x,v_{1},v_{2},\cdots,v_{s+1}\right) $ over
$F$ since $\left( z\right) $ is a (generic) specialization of
$\left( x\right) $ over $F$ and
$\left(v_{1},v_{2},\cdots,v_{s+1}\right) $ is a (generic)
specialization of $\left( v_{1},v_{2},\cdots,v_{s+1}\right) $ itself
over $F$.

\emph{Case (ii)}. Let $v_{s+1}\in \overline{F}$.

By the above assumption for $z$ it is seen that $z$ is not contained
in the field $F(v_{s+1})$. It is easily seen that there is an
isomorphism between the fields $F(x,v_{s+1})$ and $F(z,v_{s+1})$. It
follows that $(z,v_{s+1})$ is a (generic) specialization of
$(x,v_{s+1})$ over $F$. By the same procedure as in the above
\emph{Case (i)} it is seen that $\left(
z,v_{s+1},v_{1},v_{2},\cdots,v_{s}\right) $ is a (generic)
specialization of $\left( x,v_{s+1},v_{1},v_{2},\cdots,v_{s}\right)
$ over $F$.

Now take any polynomial $f\left( X,X_{1},X_{2}
,\cdots,X_{s+1}\right)$ over $F$. According to \emph{Cases
(i)-(ii)}, it is seen that $f\left(
x,v_{1},v_{2},\cdots,v_{s+1}\right) =0$ holds if and only if
$f\left( z,v_{1},v_{2},\cdots,v_{s+1}\right) =0$ holds by the theory
for generic specializations. This completes the proof.
\end{proof}

\textbf{Claim$^{\dag\dag}$}. Assume that $F(u)$ and $F(u^{\prime})$
are isomorphic over $F$ given by $u\mapsto u^{\prime}$. Let $w$ be
an element contained in an extension of $F$. Then there is an
element $w^{\prime}$ contained in some extension of $F$ such that
the fields $F(u,w)$ and $F(u^{\prime},w^{\prime})$ are isomorphic
over $F$.

\begin{proof}
It is immediate from \emph{Proposition 4} in \cite{weil}, \emph{Page 30}.
\end{proof}

\begin{corollary}
Let $K$ be a finitely generated extension of a field $k$. Then $K$
is a quasi-galois field over $k$ if and only if $K$ has one and only
one conjugation over $k$.
\end{corollary}

\begin{proof}
Prove $\Leftarrow$. Let $K$ have only one $k-$conjugation $H$. We
must have $H=K$ and then each $k-$conjugation of $K$ is contained in
$K$. By \emph{Theorem 3.12} it is seen that $K$ is a quasi-galois
field over $k$.

Prove $\Rightarrow$. Let $K$ be a $k-$quasi-galois field and $H$ a
$k-$conjugation of $K$. Choose a $k-$isomorphism $\tau$ of $H$ onto
$K$ and a $(s,m)-$nice $k-$basis $v_{1},v_{2},\cdots ,v_{m}$ of $K$
such that $H$ is a conjugate of $K$ over $F$ by $\tau$, where
$F\triangleq k(v_{1},v_{2},\cdots ,v_{s})$. We have $F\subseteq
H\subseteq K$.

Hypothesize $H\subsetneqq K$. Fixed any $x_{0}\in K\setminus H$. For
the element $x_{0}$ there are two cases.

\emph{Case (i)}. Let $x_{0}$ be a variable over $H$. We have
$$\dim_{k}H=\dim_{k}K=s< \infty$$ since $H$ and $K$ are conjugations
over $k$. But from $x_{0}\in K\setminus H$, it is seen that
$$1+\dim_{k}H=\dim_{k}H(x_{0})\leq\dim_{k}K$$ hold; from it we will obtain a
contradiction.

\emph{Case (ii)}. Let $x_{0}$ be algebraic over $H$. As $\overline
{H}\subseteq \overline{F}$, we have $x_{0}\in \overline{F}$; it
follows that $x_{0}$ is algebraic over $F$. It is clear that we have
$$[H:F]=[K:F]< \infty$$ since $H$ is a conjugate of $K$ over $F$ by
$\tau$. But from $x_{0}\in K\setminus H$, it is seen that
$$2+[H:F]\leq[H(x_{0}):F]\leq[K:F]$$ hold; from it we will obtain a
contradiction.

Therefore, the set $K\setminus H$ is empty and we must have $K=H$.
\end{proof}

\subsection{Definition for conjugations of an open set}

The notion on conjugations of an open set in a given variety that
will be defined in this subsection can be regarded as a geometric
counterpart to that for the case of fields in \S 3.3.

Let us first consider the case for integral domains. Here we let
$Fr\left( D\right) $ denote the fractional field of an integral
domain $D$.

\begin{definition}
Let $D\subseteq D_{1}\cap D_{2}$ be three integral domains.

$(i)$ The ring $D_{1}$ is said to be \textbf{$D-$quasi-galois} (or,
\textbf{quasi-galois} over $D$) if the field $Fr\left(D_{1}\right) $
is a quasi-galois extension of $Fr\left( D\right)$.

$(ii)$ Assume that there is a $(r,n)-$nice $k-$basis
$w_{1},w_{2},\cdots ,w_{n}$ of the field $Fr(D_{1})$ and an
$F-$isomorphism $\tau_{(r,n)}:Fr(D_{1})\rightarrow Fr(D_{2})$ of
fields such that  $\tau_{(r,n)}(D_{1})=D_{2}$, where $k=Fr(D)$ and
we set $F\triangleq k(w_{1},w_{2},\cdots ,w_{r})$ to be contained in
the intersection $Fr(D_{1})\cap Fr(D_{2})$.

Then the ring $D_{1}$ is said to be a \textbf{$D-$conjugation} of
the ring $D_{2}$ (or, a \textbf{conjugation} of $D_{2}$ over $D$).
\end{definition}

Then consider an integral scheme $Z$. Let $z\in Z$. By the structure
sheaf $\mathcal{O}_{Z}$ on $Z$, we have the canonical embeddings
$$i^{Z}_{U}:\mathcal{O}_{Z}(U)\rightarrow k(Z);$$
$$i^{Z}_{z}:\mathcal{O}_{Z,z}\rightarrow k(Z);$$
$$i^{z}_{U}:\mathcal{O}_{Z}(U)\rightarrow \mathcal{O}_{Z,z}$$
for every open set $U$ of $Z$ containing $z$, where
$k\left(Z\right)=\mathcal{O}_{Z,\xi }$ is the function field of $X$
and $\xi$ is the generic point of $Z$.

We will identify these integral domains with their images, that is,
we will take the rings
\begin{equation*}
\mathcal{O}_{Z}(U)\subseteq \mathcal{O}_{Z,z}\subseteq
k\left(Z\right)
\end{equation*}
as subrings of the function field $k\left(Z\right)$. This leads us
to obtain the following definitions.

Now fixed any two $k-$varieties (or, arithmetic varieties) $X$ and
$Y$ and let $\phi :X\rightarrow Y$ be a morphism of finite type.
Take a point $y \in \phi(X)$ and an open set $V$ in $Y$ with $V \cap
\phi(X) \neq \emptyset$.

\begin{definition}
Assume that $U_{1}$ and $U_{2}$ are open sets of $X$ such that
either $U_{1}$ or $U_{2}$ is contained in $\phi ^{-1}(V)$. The open
set $U_{1}$ is said to be a \textbf{$V-$conjugation} of the open set
$U_{2}$ if the ring $i^{X}_{U_{1}}(\mathcal{O}_{X}(U_{1}))$
($\subseteq k(X)$) is a conjugation of the ring
$i^{X}_{U_{2}}(\mathcal{O}_{X}(U_{2}))$ ($\subseteq k(X)$) over the
ring $i^{X}_{\phi^{-1}(V)}(\phi^{\sharp}(\mathcal{O}_{Y}(V)))$
($\subseteq k(X)$), where
$\phi^{\sharp}:\mathcal{O}_{Y}(V)\rightarrow
\phi_{\ast}\mathcal{O}_{X}(V)=\mathcal{O}_{X}(\phi^{-1}(V))$ is the
ring homomorphism.

If $U_{1}$ and $U_{2}$ are both contained in $\phi ^{-1}(V)$, such a
$V-$conjugation is said to be \textbf{geometric}.
\end{definition}

\begin{remark}
It is seen that the above conjugation of an open set is well-defined
since we have
\begin{equation*}
\begin{array}{l}
i^{X}_{\phi^{-1}(V)}(\phi^{\sharp}(\mathcal{O}_{Y}(V)))\\

=i^{X}_{U_{1}}(i_{\phi^{-1}(V)}^{U_{1}}(\phi^{\sharp}(\mathcal{O}_{Y}(V))))\\

=i^{X}_{U_{1}\cap U_{2}}(i_{\phi^{-1}(V)}^{U_{1}\cap U_{2}}(\phi^{\sharp}(\mathcal{O}_{Y}(V))))\\

=i^{X}_{U_{2}}(i_{\phi^{-1}(V)}^{U_{2}}(\phi^{\sharp}(\mathcal{O}_{Y}(V)))).
\end{array}
\end{equation*}
In particular, if $\phi$ is surjective, we have
$$\phi^{\sharp}(k(Y))\subseteq k(X);$$
$$\phi^{\sharp}(i^{Y}_{V} (\mathcal{O}_{Y}(V))
=i^{X}_{\phi^{-1}(V)}(\phi^{\sharp}(\mathcal{O}_{Y}(V))).$$
\end{remark}

\begin{definition}
Assume that either $x_{1}\in X$ or  $x_{2}\in X$ is contained in $
\phi ^{-1}\left( y\right) $. The point $x_{1}$ is said to be a
\textbf{$y-$conjugation} of the point $x_{2}$ if the ring
$i^{X}_{x_{1}}\left( \mathcal{O}_{X,x_{1}}\right) $ ($\subseteq
k(X)$) is a conjugation of the ring $i^{X}_{x_{2}}\left(
\mathcal{O}_{X,x_{2}}\right) $ ($\subseteq k(X)$) over the ring
$i^{X}_{x_{1}}(\phi^{\sharp}(\mathcal{O}_{Y,y}))$ ($\subseteq
k(X)$), where $\phi^{\sharp}:\mathcal{O}_{Y,y}\rightarrow
\mathcal{O}_{X,x_{1}}$ is the ring homomorphism.

If $x_{1}$ and $x_{2}$ are both contained in $\phi ^{-1}(y)$, such a
$y-$conjugation is said to be \textbf{geometric}.
\end{definition}

\begin{remark}
The above conjugation of a point is well-defined. In deed, by
\emph{Remark 3.17} we have
$$i^{X}_{x_{1}}(\phi^{\sharp}(\mathcal{O}_{Y,y}))
=i^{X}_{x_{2}}(\phi^{\sharp}(\mathcal{O}_{Y,y}))$$ as subrings of $k(X)$
according to
the preliminary facts on direct systems of rings.
\end{remark}

Let $A$ be a commutative ring with identity. $A$ is said to be
\textbf{affinely realized} in $X$ by an open set $U$ of $X$ if we
have $A=\mathcal{O}_{X}(U)$. $A$ is said to be \textbf{affinely
realized} in $X$ by a point $x$ of $X$ if we have $A=
\mathcal{O}_{X,x}$. This is a hint of the following notion for the
case of varieties.

\begin{definition}
An open set $U\subseteq \phi ^{-1}(V)$ in the variety $X$ is said to
have a \textbf{quasi-galois set of $V-$conjugations} in $X$ if each
conjugation $A$ of the ring $i^{X}_{U}(\mathcal{O}_{X}(U))$ over the
ring $i^{X}_{\phi^{-1}(V)}(\phi^{\sharp}(\mathcal{O}_{Y}(V)))$ can
be affinely realized canonically by an open set $U_{A}$ of $X$ such
that $A=i^{X}_{U_{A}}(\mathcal{O}_{X}(U_{A})).$
\end{definition}

It is easily seen that such an open set $U_{A}$ can be contained in
the set $\phi ^{-1}(V)$.

\begin{definition}
A point $x\in \phi ^{-1}\left( y\right) $ in the variety $X$ is said
to have a \textbf{quasi-galois set of $y-$conjugations} in $X$ if
each conjugation $A$ of the ring $i^{X}_{x}\left(
\mathcal{O}_{X,x}\right) $ over the ring
$i^{X}_{x}(\phi^{\sharp}(\mathcal{O}_{Y,y}))$ can be affinely
realized canonically by a point $x_{A}$ of $X$ such that
$A=i^{X}_{x_{A}}\left( \mathcal{O}_{X,x_{A}}\right) $.

In particular, the fiber $\phi ^{-1}\left( y\right)$ is said to be
\textbf{quasi-galois} over $y$ if each point of the fiber $\phi
^{-1}\left( y\right)$ has a quasi-galois set of $y-$conjugations in
$X$.
\end{definition}

\begin{remark}
Let $y\in V$. By \emph{Theorem 3.23} below it is easily seen that
each point $x_0\in \phi ^{-1}\left( y\right) $ has a quasi-galois
set of $y-$conjugations implies that each affine open set
$U\subseteq \phi ^{-1}(V)$ containing $x_0$ has a quasi-galois set
of $V-$conjugations in $X$.
\end{remark}

\subsection{Quasi-galois closed varieties and conjugations of open sets}

In this subsection we will obtain some properties of quasi-galois
closed varieties by virtue of conjugations of open sets.

\begin{theorem}
Let $X$ and $Y$ be two $k-$varieties (or, two arithmetic varieties)
such that $X$ is quasi-galois closed over $Y$ by
 a surjective morphism $\phi$ of finite
type.

$(i)$ Fixed any affine open set $V$ of $Y$. Then each affine open
set $U\subseteq \phi ^{-1}(V)$ has a quasi-galois set of
$V-$conjugations in $X$.

$(ii)$ Let $\Omega$ be an fixed algebraic closure of the functional
field $k(X)$. Then  we have
 $$\mathcal{O}_{X}(U)\subseteq \mathcal{O}_{X,x_{0}}
 \subseteq \Omega$$ exactly as subsets for any point $x\in X$
and any affine open set $U$ of $X$ containing $x$.
\end{theorem}

\begin{proof}
Let $\Omega$ be an fixed algebraic closure of the functional field
$k(X)$ of $X$. By \emph{Proposition 3.9} we have the natural affine
structure $\mathcal{A}(\mathcal{O}_{X})$ of the variety
$(X,\mathcal{O}_{X})$ such that $\mathcal{A}(\mathcal{O}_{X})$ is
with values in $\Omega$ and that $\mathcal{O}_{X}$ is an extension
of $\mathcal{A}(\mathcal{O}_{X})$.

$(i)$ Hypothesize that there is an affine open set $U_{0}\subseteq
\phi ^{-1}(V)$ such that a conjugation $H$ of
$i^{X}_{U_{0}}(\mathcal{O}_{X}(U_{0}))$  over
$i^{X}_{\phi^{-1}(V)}(\phi^{\sharp}(\mathcal{O}_{Y}(V)))$ can not be
affinely realized canonically by any open set $U^{\prime}$ of $X$
with $H=i^{X}_{U^{\prime}}(\mathcal{O}_{X}(U^{\prime})).$

Evidently, $H\not= i^{X}_{U_{0}}(\mathcal{O}_{X}(U_{0}))$. From the
field $\Omega$ we have
$\mathcal{O}_{X}(U_{0})=i^{X}_{U_{0}}(\mathcal{O}_{X}(U_{0}))$ and
then $H\not= \mathcal{O}_{X}(U_{0})$.

Put
$$\mathcal{C}^{\prime}_{X}=\{(U_{0},\phi^{\prime}_{0};H)\}
\bigcup(\mathcal{A}(\mathcal{O}_{X}) \setminus
\{(U_{0},\phi_{0};A_{0})\})$$ where $(U_{0},\phi_{0};A_{0})\in
\mathcal{A}(\mathcal{O}_{X})$ and $\phi^{\prime}_{0}(U_{0})=Spec(H)$
is an isomorphism.

Let $\Gamma(\mathcal{C}^{\prime}_{X})$ be the maximal pseudogroup of
affine transformations in
$(X,\mathcal{O}_{X};\mathcal{C}^{\prime}_{X})$ and let
$\mathcal{A}^{\prime}(\mathcal{O}_{X})$ be the affine
$\Gamma(\mathcal{C}^{\prime}_{X})-$structure defined by the reduced
affine covering $\mathcal{C}^{\prime}_{X}$.

By gluing schemes (see \cite{Hrtsh}), it is easily seen that
$\mathcal{A}^{\prime}(\mathcal{O}_{X})$ is admissible and there is a
sheaf $\mathcal{O}^{\prime}_{X}$ on $X$ such that
$\mathcal{O}^{\prime}_{X}$ is an extension of
$\mathcal{A}^{\prime}(\mathcal{O}_{X})$. Then
$(X,\mathcal{O}^{\prime}_{X})$ is a scheme such that
$\mathcal{O}^{\prime}_{X}(U_{0})$ is exactly equal to the ring $H$
since they are both subrings of $\Omega$.

As $\mathcal{A}(\mathcal{O}_{X}) $ and
$\mathcal{A}^{\prime}(\mathcal{O}_{X})$ are both  with values in
$\Omega$, in virtue of $(ii)$ of \emph{Proposition 3.9} we have
$\mathcal{A}(\mathcal{O}_{X})\supseteq
\mathcal{A}^{\prime}(\mathcal{O}_{X}) $; as affine structures are
reduced coverings of $X$, we must have
$$\mathcal{O}_{X}(U_{0})=\mathcal{O}^{\prime}_{X}(U_{0})=H$$ since
$(U_{0},\phi^{\prime}_{0};H)\in\mathcal{A}^{\prime}(\mathcal{O}_{X})$,
which will be in contradiction with the hypothesis above. Therefore,
each affine open set $U\subseteq \phi ^{-1}(V)$ has a quasi-galois
set of $V-$conjugations in $X$.

$(ii)$ Fixed a point $x_{0}\in X$. Let $I_{x_0}$ (respectively,
$J_{x_0}$) be the index family of open sets (respectively, affine
open sets) of $X$ containing $x_0$. By set inclusion, $I_{x_0}$ and
$J_{x_0}$ are partially ordered sets and then are directed sets. It
is easily seen that $J_{x_0}$ and $I_{x_0}$ are cofinal since affine
open sets for a base for the topology on the space of $X$. Hence,
the stalk $\mathcal{O}_{X,x_0}$ at $x_0$ is the direct limit of the
system of rings $\mathcal{O}_{X}(U)$ with $U\in J_{x_0}$.

Now consider the natural affine structure
$\mathcal{A}(\mathcal{O}_{X})$ with values in $\Omega$. Take each
local chart $(U,\phi;A)\in \mathcal{A}(\mathcal{O}_{X})$ with $U \in
J_{x_0}$. It is seen that we have
$$i^{X}_{U}(\mathcal{O}_{X}(U))=\mathcal{O}_{X}(U)=A \subseteq
\Omega$$ since $\mathcal{O}_{X}$ is an extension of
$\mathcal{A}(\mathcal{O}_{X})$.

Similarly, take any affine open sets $U_{1},U_{2}$ of $X$ containing
$x_{0}$ such that $ U_{1}\subseteq U_{2}$. We have
$$i^{U_{2}}_{U_{1}}(\mathcal{O}_{X}(U_{2}))=\mathcal{O}_{X}(U_{2})
\subseteq \Omega.$$

 It follows that for the stalk of $\mathcal{O}_{X}$ at $x_0$ we have
 $$\mathcal{O}_{X,x_{0}}=\bigcup_{U \in
J_{x_0}}\mathcal{O}_{X}(U)\subseteq \Omega.$$

(Please notice that all isomorphisms $i^{X}_{U}$ and
$i^{U_{2}}_{U_{1}}$ here are exactly identity maps only for affine
open sets!)
\end{proof}

\begin{theorem}
Let $X$ and $Y$ be two $k-$varieties (or, two arithmetic varieties)
such that $X$ is quasi-galois closed over $Y$ by
 a surjective morphism $\phi$ of finite
type. Then the function field $k\left( X\right) $ is a quasi-galois
extension over the image $\phi^{\sharp}(k\left( Y\right)) $ of the
function field $f(Y)$.
\end{theorem}

\begin{proof}
By  \emph{Proposition 3.9} it is seen that there is the natural
affine structure $\mathcal{A}(\mathcal{O}_{X})$ of the variety
$(X,\mathcal{O}_{X})$ with values in $\Omega$ and $\mathcal{O}_{X}$
is an extension of $\mathcal{A}(\mathcal{O}_{X})$, where $\Omega$ is
an fixed algebraic closure of the functional field $k(X)$ of $X$.

Now fixed any conjugation $H$ of $ k\left( X\right) $ over
$\phi^{\sharp}(k(Y))$. Take any element $w_{0}\in H$. Let $\sigma
:H\rightarrow k\left( X\right) $ be an isomorphism over
$\phi^{\sharp}(k(Y))$. Put $ u_{0}=\sigma \left( w_{0}\right) . $

In virtue of \emph{Theorem 3.23} we have
 $$\mathcal{O}_{X}(U)\subseteq k(X)=\mathcal{O}_{X,\xi}
 \subseteq \Omega$$ exactly as subsets of $\Omega$ for
  any affine open set $U$ of $X$, where $\xi$ is
the generic point of $X$. Then we have
$$\bigcup_{U}\mathcal{O}_{X}(U)=\mathcal{O}_{X,\xi}$$ since affine
open sets $U$ of $X$ form a base for the topology of $X$.

It follows that there are affine open subsets $V_0$ of $Y$ and
$U_0\subseteq \phi ^{-1}\left( V_)\right) $ of $X$ such that $u_{0}$
is contained in $\mathcal{O}_{X}\left( U_0\right) $. By
\emph{Theorem 3.23} again it is seen that there is some affine open
set $W_0$ of $X$ such that $W_0$ is a $V-$conjugation of $U_0$ and
that the element $w_0$ is contained in $\mathcal{O}_{X}(W_0)$.

Hence, $w_0$ is contained in
$k(X)=\mathcal{O}_{X,\xi}\subseteq\Omega$. This proves any given
conjugation $H$ of $ k\left( X\right) $ over $\phi^{\sharp}(k(Y))$
is contained in $k(X)$. From \emph{Theorem 3.12} it is seen that the
function field $k(X)$ is a quasi-galois extension of the field
$\phi^{\sharp}(k(Y))$.
\end{proof}

At last we have the following corollary.

\begin{corollary}
Let $X$ and $Y$ be two $k-$varieties (or, two arithmetic varieties)
such that $X$ is quasi-galois closed over $Y$ by
 a surjective morphism $\phi$ of finite
type. Suppose that each $V-$conjugation of $U$ is geometric for any
affine open sets $V \subseteq Y$ and $U\subseteq \phi ^{-1} (V)$.

Then each point $x_0\in \phi ^{-1}(y_0)$ has a quasi-galois set of
geometric $y-$conjugations in $X$ for any point $y_0\in Y$.
\end{corollary}

\begin{proof}
Fixed a  point $y_0\in Y$ and a point $x_0\in \phi ^{-1}(y_0)$. Take
any affine open sets $V \subseteq Y$ and $U\subseteq \phi ^{-1} (V)$
such that $x_0\in U$ and $y_0\in V$. By \emph{Theorem 3.23} it is
seen that $U$ has a quasi-galois set of geometric $V-$conjugations;
from \emph{Theorem 3.24} it is seen that each point $x_0\in \phi
^{-1}(y_0)$ has a quasi-galois set of geometric $y-$conjugations in
$X$ for any point $y_0\in Y$.
\end{proof}

\subsection{Automorphism groups of quasi-galois closed varieties and Galois groups of the function fields}

For automorphism groups of  quasi-galois closed varieties, we have
the following result.

\begin{theorem}
Let $X$ and $Y$ be two $k-$varieties (or, two arithmetic varieties)
such that $X$ is quasi-galois closed over $Y$ by
 a surjective morphism $\phi$ of finite
type. Suppose that $k\left( X\right) /\phi^{\sharp}(k(Y))$ is
separably generated. Then the function field $k\left( X\right) $ is
a Galois extension of $\phi^{\sharp}(k(Y))$ and there is a group
isomorphism
$$
{Aut}\left( X/Y\right) \cong Gal(k\left( X\right) /\phi^{\sharp}(k(Y))).
$$
\end{theorem}

\begin{proof}
$\textbf{(i).}$ Prove that the function field $k\left( X\right)$ is
a finitely generated Galois extension of $\phi^{\sharp}(k(Y))$.

Without loss of generality, assume that $k\left( X\right)$ is a
transcendental extension over $\phi^{\sharp}(k(Y))$. By
\emph{Corollary 3.14} and \emph{Theorem 3.24} it is seen that every
conjugation of $k\left( X\right)$ over $\phi^{\sharp}(k(Y))$ is
$k\left( X\right)$ itself. It needs to prove that there exists a
$\sigma_{0} \in Gal(k\left( X\right)/\phi^{\sharp}(k(Y)))$ such that
$\phi^{\sharp}(k(Y))$ is the invariant subfield of $\sigma_{0}$.

In deed, take a $(r,n)-$nice $F-$basis $v_{1},v_{2},\cdots ,v_{n}$
of $k\left( X\right)$. By the assumption above it is seen that
$r\geqslant 1$ holds and $k\left( X\right)$ is an algebraic Galois
extension of the field $F_{0}\triangleq
\phi^{\sharp}(k(Y))(v_{1},v_{2},\cdots ,v_{r})$. Fixed any
$\tau_{0}\in Gal(k\left( X\right)/F_{0})$ with $\tau_{0}\not=
id_{k\left( X\right)}$. Let $\tau_{1}\in
Gal(F_{0}/\phi^{\sharp}(k(Y)))$ be given by

$$v_{1}\mapsto \frac{1}{v_{1}},
v_{2}\mapsto \frac{1}{v_{2}},\cdots,v_{r}\mapsto \frac{1}{v_{r}}.$$

We have a $\sigma_{0} \in Gal(k\left( X\right)/\phi^{\sharp}(k(Y)))$
defined by $\tau_{0}$ and $\tau_{1}$ in such a manner
$$\frac{f(v_{1},v_{2},\cdots ,v_{n})}{g(v_{1},v_{2},\cdots
,v_{n})}\in k(X)$$
$$\mapsto \frac{f(\tau_{1}(v_{1}),\tau_{1}(v_{2}),\cdots
,\tau_{1}(v_{r}),\tau_{0}(v_{r+1}),\cdots,
\tau_{0}(v_{n}))}{g(\tau_{1}(v_{1}),\tau_{1}(v_{2}),\cdots
,\tau_{1}(v_{r}),\tau_{0}(v_{r+1}),\cdots, \tau_{0}(v_{n}))}\in k(X)$$ for
any polynomials $f(X_{1},X_{2},\cdots, X_{n})$ and
$g(X_{1},X_{2},\cdots, X_{n})\not= 0$ over the field
$\phi^{\sharp}(k(Y))$ such that $g(v_{1},v_{2},\cdots ,v_{n})\not=
0$. By $\tau_{0}$ it is seen that we have $$g(v_{1},v_{2},\cdots
,v_{n})= 0$$ if and only if
$$g(\tau_{1}(v_{1}),\tau_{1}(v_{2}),\cdots
,\tau_{1}(v_{r}),\tau_{0}(v_{r+1}),\cdots, \tau_{0}(v_{n}))=0$$
holds. Hence, $\sigma_{0}$ is well-defined.

It is seen that $\phi^{\sharp}(k(Y))$ is the invariant subfield of
$\sigma_{0}$ and then $\phi^{\sharp}(k(Y))$ is the invariant
subfield of $Gal(k\left( X\right)/\phi^{\sharp}(k(Y)))$. Therefore,
$k\left( X\right)$ is a Galois extension of $\phi^{\sharp}(k(Y))$.

$\textbf{(ii).}$ Now let $\mathcal{A}(\mathcal{O}_{X})$ be the
natural affine structure of the variety $(X,\mathcal{O}_{X})$ with
values in
 an fixed algebraic
closure $\Omega$ of the functional field $k(X)$.

For an open set $H$ in $X$, we have an isomorphism
$\tau_{H}:\Gamma(H,\mathcal{O}_{X}) \cong \mathcal{O}_{X}(H)$ of
algebras and an embedding $i^{X}_{H}:\mathcal{O}_{X}(H) \rightarrow
\mathcal{O}_{X,\xi}\subseteq \Omega$, where $\xi$ is the generic
point of $x$.

For the sake of convenience,  all such rings
$\Gamma(H,\mathcal{O}_{X})$ and $\mathcal{O}_{X}(H)$ are regarded as
the subrings of the function field $k(X)$ by the maps
$i^{X}_{H}\circ \tau_{H}$ and $i^{X}_{H}$, respectively.

The function field $k(X)=\mathcal{O}_{X,\xi}$  is regarded as the
set of the elements of the forms $ (U,f) $, where $U$ is an open set
of $X$ and $f$ is an element of $\mathcal{O}_{X}(U)$ ($\subseteq
\Omega$). That is,
$$k(X)=\{(U,f):f\in \mathcal{O}_{X}(U)\text{ and }U\subseteq X\text{
is open}\}.$$

In the following we will proceed in several steps to demonstrate
that there exists an isomorphism
\begin{equation*}
t:{Aut}\left( X/Y\right) \cong {Gal}\left( k\left( X\right)
/\phi^{\sharp}\left(k\left( Y\right)\right) \right)
\end{equation*}
of groups.

\emph{Step 1.} Fixed any automorphism $ \sigma =\left( \sigma
,\sigma ^{\sharp}\right) \in Aut \left( X/Y\right) . $ That is, $
\sigma : X \longrightarrow X $ is a homeomorphism and $ \sigma ^
{\sharp}:\mathcal{O}_{X} \rightarrow \sigma _{\ast }\mathcal{O}_{X}
$ is an isomorphism of sheaves of rings on $X$. As $\dim X<\infty$,
we have $\sigma (\xi)=\xi$. It follows that
\begin{equation*}
\sigma ^{\sharp}:k\left( X\right)=\mathcal{O}_{X,\xi } \rightarrow \sigma
_{\ast }\mathcal{O}_{X,\xi }=k\left( X\right)
\end{equation*}
is an automorphism of $k(X)$. Let $\sigma ^{\sharp-1}$ denote the inverse of $%
\sigma ^{\sharp}$.

Take any open subset $U$ of $X$. We have the restriction
\begin{equation*}
\sigma=(\sigma ,\sigma ^{\sharp}): (U,\mathcal{O}_{X}|_{U}) \longrightarrow
(\sigma(U),\mathcal{O}_{X}|_{\sigma(U)})
\end{equation*}
of open subschemes. That is,
\begin{equation*}
\sigma^{\sharp}:\mathcal{O}_{X}|_{\sigma(U)} \rightarrow \sigma_{\ast}%
\mathcal{O}_{X}|_{U}
\end{equation*}
is an isomorphism of sheaves on $\sigma (U)$. In particular,
\begin{equation*}
\sigma^{\sharp}:\mathcal{O}_{X}(\sigma(U))
=\mathcal{O}_{X}|_{\sigma(U)}(\sigma(U)) \rightarrow \mathcal{O}_{X}(U)=\sigma_{\ast}
\mathcal{O}_{X}|_{U}(\sigma(U))
\end{equation*}
is an isomorphism of rings.

For every $f \in \mathcal{O}_{X}|_{U}(U)$, we have
\begin{equation*}
f \in \sigma_{\ast}\mathcal{O}_{X}|_{U}(\sigma(U));
\end{equation*}
hence
\begin{equation*}
\sigma^{\sharp -1}(f) \in \mathcal{O}_{X}(\sigma(U)).
\end{equation*}

Now define a mapping
\begin{equation*}
t:Aut \left( X/Y\right) \longrightarrow Gal\left( k\left( X\right)
/\phi^{\sharp}(k\left( Y\right) )\right)
\end{equation*}
given by
\begin{equation*}
\sigma =(\sigma ,\sigma ^{\sharp})\longmapsto t(\sigma)=\left\langle \sigma ,\sigma ^{\sharp
-1}\right\rangle
\end{equation*}
such that
\begin{equation*}
\left\langle \sigma ,\sigma ^{\sharp-1}\right\rangle :\left( U,f\right)\in \mathcal{O}_{X}(U)
\longmapsto \left( \sigma \left( U\right) ,\sigma ^{\sharp-1}\left( f\right)
\right)\in \mathcal{O}_{X}(\sigma(U))
\end{equation*}
is a mapping of $k(X)$ into $k(X)$.

\emph{Step 2.} Prove that $t$ is well-defined. In deed, given any
\begin{equation*}
\sigma =\left( \sigma ,\sigma ^{\sharp}\right) \in Aut\left( X/Y\right).
\end{equation*}
For any $(U,f),(V,g) \in k(X)$, we have
\begin{equation*}
(U,f)+(V,g)=(U\cap V, f+g)
\end{equation*}
and
\begin{equation*}
(U,f)\cdot(V,g)=(U\cap V, f\cdot g);
\end{equation*}
then we have
\begin{equation*}
\begin{array}{l}
\left\langle \sigma ,\sigma ^{\sharp-1}\right\rangle((U,f)+(V,g)) \\
=\left\langle \sigma ,\sigma ^{\sharp-1}\right\rangle((U\cap V, f+g)) \\
=(\sigma(U\cap V), \sigma^{\sharp -1}(f+g)) \\
=(\sigma(U\cap V), \sigma^{\sharp -1}(f))+(\sigma(U\cap V), \sigma^{\sharp
-1}(g)) \\
=(\sigma(U), \sigma^{\sharp -1}(f))+(\sigma(V), \sigma^{\sharp -1}(g)) \\
=\left\langle \sigma ,\sigma ^{\sharp-1}\right\rangle((U,f))+ \left\langle
\sigma ,\sigma ^{\sharp-1}\right\rangle((V,g))
\end{array}%
\end{equation*}
and
\begin{equation*}
\begin{array}{l}
\left\langle \sigma ,\sigma ^{\sharp-1}\right\rangle((U,f)\cdot(V,g)) \\
=\left\langle \sigma ,\sigma ^{\sharp-1}\right\rangle((U\cap V, f\cdot g))
\\
=(\sigma(U\cap V), \sigma^{\sharp -1}(f\cdot g)) \\
=(\sigma(U\cap V), \sigma^{\sharp -1}(f))\cdot(\sigma(U\cap V),
\sigma^{\sharp -1}(g)) \\
=(\sigma(U), \sigma^{\sharp -1}(f))\cdot(\sigma(V), \sigma^{\sharp -1}(g))
\\
=\left\langle \sigma ,\sigma ^{\sharp-1}\right\rangle((U,f))\cdot
\left\langle \sigma ,\sigma ^{\sharp-1}\right\rangle((V,g)).%
\end{array}%
\end{equation*}

It follows that $\left\langle \sigma ,\sigma ^{\sharp-1}\right\rangle$ is an
automorphism of $k\left( X\right) .$

It needs to prove that $\left\langle \sigma ,\sigma
^{\sharp-1}\right\rangle$ is an isomorphism over
$\phi^{\sharp}(k(Y))$. In deed, consider the given morphism
$\phi=(\phi,\phi^{\sharp}):(X,\mathcal{O} _{X})\rightarrow
(Y,\mathcal{O}_{Y})$ of schemes. Evidently, $\phi(\xi)$ is the
generic point of $Y$ and $\xi$ is invariant under any automorphism
$\sigma \in Aut\left( X/Y\right)$. Then $\sigma^{\sharp}:
\mathcal{O}_{X,\xi} \rightarrow \mathcal{O}_{X,\xi}$ is an
isomorphism of algebras over
$\phi^{\sharp}(\mathcal{O}_{Y,\phi(\xi)})=\phi^{\sharp}(k(Y))$.
Hence, $\left\langle \sigma ,\sigma
^{\sharp-1}\right\rangle|_{\phi^{\sharp}(k(Y))}=id_{\phi^{\sharp}(k(Y))}.$

This proves
$$
\left\langle \sigma ,\sigma ^{\sharp-1}\right\rangle \in Gal\left( k\left(
X\right) /\phi^{\sharp}(k\left( Y\right) \right)) .
$$ That is, $t$ is a well-defined map.

Prove that $t$ is a homomorphism between groups. In fact, take any
\begin{equation*}
\sigma =\left( \sigma ,\sigma ^{\sharp}\right) ,\delta =\left( \delta
,\delta ^{\sharp}\right) \in Aut\left( X/Y\right) .
\end{equation*}
By preliminary facts on schemes (see \cite{EGA}) we have
\begin{equation*}
\delta ^{\sharp -1}\circ \sigma^{\sharp -1}=(\delta \circ \sigma)^{\sharp -1};
\end{equation*}
then
\begin{equation*}
\left\langle \delta ,\delta ^{\sharp-1}\right\rangle \circ \left\langle
\sigma ,\sigma ^{\sharp-1}\right\rangle =\left\langle \delta \circ \sigma
,\delta ^{\sharp-1}\circ \sigma ^{\sharp-1}\right\rangle.
\end{equation*}

Hence, the map
\begin{equation*}
t:Aut\left( X/Y\right) \rightarrow Gal\left( k\left( X\right)
/\phi^{\sharp}\left(k\left( Y\right)\right) \right)
\end{equation*}
is a homomorphism of groups.

\emph{Step 3.} Prove that ${t}$ is injective. Assume $\sigma ,\sigma
^{\prime }\in {Aut}\left( X/Y\right) $ such that $t\left( \sigma
\right) =t\left( \sigma ^{\prime }\right) .$ We have
\begin{equation*}
\left( \sigma \left( U\right) ,\sigma ^{\sharp-1}\left( f\right) \right)
=\left( \sigma ^{\prime }\left( U\right) ,\sigma ^{\prime \sharp-1}\left(
f\right) \right)
\end{equation*}
for any $\left( U,f\right) \in k\left( X\right) .$ In particular, we
have
\begin{equation*}
\left( \sigma \left( U_{0}\right) ,\sigma ^{\sharp-1}\left( f\right) \right)
=\left( \sigma ^{\prime }\left( U_{0}\right) ,\sigma ^{\prime
\sharp-1}\left( f\right) \right)
\end{equation*}
for any $f\in \mathcal{O}_{X}(U_{0})$ and any affine open subset
$U_{0}$ of $X$ such that $ \sigma \left( U_{0}\right)$ and $\sigma
^{\prime }\left( U_{0}\right)$ are both contained in $\sigma \left(
U\right) \cap \sigma ^{\prime }\left( U\right) $.

As $\mathcal{O}_{X}$ is an extension of
$\mathcal{A}(\mathcal{O}_{X})$, there are three subrings
$$A_0=\mathcal{O}_{X}(U_{0}),B_0=\mathcal{O}_{X}(\sigma(U_{0})),
\text{ and } {B_{0}}^{\prime}=\mathcal{O}_{X}(\sigma^{\prime}(U_{0}))$$ of $k(X)$ such that
$$B_0=\sigma^{\sharp-1}(A_0) =\sigma^{\prime \sharp-1}(A_0)=B_{0}^{\prime}.$$

By preliminary facts on affine schemes (see \cite{EGA}) again, it is seen
that
$$\sigma |_{U_{0}}=\sigma ^{\prime }|_{U_{0}}$$ holds as
isomorphisms of schemes. As $U_0$ is dense in $X$, we have
$$\sigma =\sigma |_{\overline{U_{0}}}=\sigma ^{\prime
}|_{\overline{U_{0}}}= \sigma ^{\prime }$$ on the whole of $X$. This proves that $t$ is
an injection.

\emph{Step 4.} Prove that ${t}$ is surjective.  Fixed any element
$\rho$ of the group $Gal\left( k\left( X\right)
/\phi^{\sharp}\left(k\left( Y\right)\right) \right) $.

As $k(X)=\{(U_{f},f):f\in \mathcal{O}_{X}(U_{f})\text{ and
}U_{f}\subseteq X\text{ is open}\}$, we have
\begin{equation*}
\rho :\left( U_{f},f\right) \in k\left( X\right) \longmapsto \left( U_{\rho
\left( f\right) },\rho \left( f\right) \right) \in k\left( X\right) ,
\end{equation*}
where $U_{f}$ and $U_{\rho (f)}$ are open sets in $X$, $f$ is
contained in $\mathcal{O}_{X}(U_{f})$, and $\rho (f)$ is contained
in $ \mathcal{O}_{X}(U_{\rho (f)})$.

We will proceed in the following several sub-steps to prove that
each element of $Gal\left( k\left( X\right)
/\phi^{\sharp}\left(k\left( Y\right)\right) \right) $  give us a
unique element of ${Aut}(X/Y)$.

$\emph{(a)}$ Fixed any affine open set $V$ of $Y$. Prove that for
each affine open set $U\subseteq \phi^{-1}(V)$ there is an affine
open set $U_{\rho}$ in $X$ such that $\rho$ determines an
isomorphism between affine schemes $(U,\mathcal{O}_{X}|_{U})$ and
$(U_{\rho},\mathcal{O}_{X}|_{U_{\rho}})$.

In fact, take any local chart $\left( U,\phi; A_{U}\right)\in
\mathcal{A}(\mathcal{O}_{X})$ with $U\subseteq \phi^{-1}(V)$ for
some affine open set $V$ of $Y$. Here $\mathcal{A}(\mathcal{O}_{X})$
is the natural affine structure of the variety $X$ with values in
$\Omega$. As $\mathcal{O}_{X}$ is an extension of
$\mathcal{A}(\mathcal{O}_{X})$, by \emph{Theorem 3.23} we have
\begin{equation*}
A= \mathcal{O}_{X}(U) =\{\left( U_{f},f\right) \in k\left( X\right)
:U_{f}\supseteq U\}
\end{equation*}
since $U$ is an affine open set of $X$. Put
\begin{equation*}
B=\{\left( U_{\rho \left( f\right) },\rho \left( f\right) \right) \in
k\left( X\right) :\left( U_{f},f\right) \in A
\}.
\end{equation*}

Then $B$ is a subring of $k(X)$. As $\rho$ is an isomorphism over
$\phi^{\sharp}(k(Y))$, it is seen that by $\rho$ the rings $A$ and
$B$ are isomorphic algebras over $\phi^{\sharp}(k(Y))$. It follows
that $A$ and $B$ are conjugations over
$\phi^{\sharp}(\mathcal{O}_{Y}(V))$.

By \emph{Theorem 3.23} again it is seen that $U$ has a quasi-galois
set of $V-$conjugations in $X$. Then there is an open set $U_{\rho}$
that is a $V-$conjugation of $U$ such that
$B=\mathcal{O}_{X}(U_{\rho})$. As $U$ is affine open, it is clear
that $U_{\rho}$ is affine open.

Hence, by $\rho$ we have a unique isomorphism
\begin{equation*}
\lambda_{U}=\left(\lambda_{U}, \lambda_{U}^{\sharp} \right): (U, \mathcal{O}%
_{X}|_{U}) \rightarrow (U_{\rho}, \mathcal{O}_{X}|_{U_{\rho}})
\end{equation*}
of the affine open subscheme in $X$ such that
\begin{equation*}
\rho |_{\mathcal{O}_{X}(U)}=\lambda_{U}^{\sharp -1}: \mathcal{O}_{X}(U)
\rightarrow \mathcal{O}_{X}(U_{\rho}).
\end{equation*}

$\emph{(b)}$ Take any affine open sets $V\subseteq Y$ and
$U,U^{\prime}\subseteq\phi^{-1}(V)$. Prove that
\begin{equation*}
\lambda_{U}|_{U\cap U^{\prime}}=\lambda_{U^{\prime}}|_{U\cap U^{\prime}}
\end{equation*}
holds as morphisms of schemes.

In fact, by the above construction for each $\lambda_{U}$ it is seen
that $\lambda_{U}^{\sharp}$ and $\lambda^{\sharp}_{U^{\prime}}$
coincide on the intersection $U\cap U^{\prime}$ since we have
\begin{equation*}
\rho |_{\mathcal{O}_{X}(U\cap U^{\prime})}=\lambda_{U}|_{U\cap U^{\prime}}^{\sharp -1}:
\mathcal{O}_{X}(U\cap U^{\prime}) \rightarrow \mathcal{O}_{X}({(U\cap U^{\prime})}_{\rho});
\end{equation*}
\begin{equation*}
\rho |_{\mathcal{O}_{X}(U\cap U^{\prime})}=\lambda_{U^{\prime}}|_{U\cap U^{\prime}}^{\sharp -1}:
\mathcal{O}_{X}(U\cap U^{\prime}) \rightarrow \mathcal{O}_{X}({(U\cap U^{\prime})}_{\rho}).
\end{equation*}

For any  point $x\in U\cap U^{\prime}$, we must have $
\lambda_{U}(x)=\lambda_{U^{\prime}}(x)$. Otherwise, if
$\lambda_{U}(x)\not=\lambda_{U^{\prime}}(x)$, will have an affine
open subset $X_0$ of $X$ that contains one of the two points
$\lambda_{U}(x)$ and $\lambda_{U^{\prime}}(x)$ but does not contain
the other since the underlying space of $X$ is a Kolmogrov space.
Assume $\lambda_{U}(x)\in X_{0}$ and $\lambda_{U^{\prime}}(x)\not
\in X_{0}$. We choose an affine open subset $U_{0}$ of $X$ such that
$x\in U_{0}\subseteq U\cap U^{\prime}$ and
$\lambda_{U}(U_{0})\subseteq X_{0}$ since we have
$$\lambda_{U}(U\cap U^{\prime})={(U\cap U^{\prime})}_{\rho}\subseteq
U_{\rho};$$
$$\lambda_{U^{\prime}}(U\cap U^{\prime})={(U\cap
U^{\prime})}_{\rho}\subseteq U^{\prime}_{\rho}.$$ However, by the
definition for each $\lambda_{U}$, we have
$$\lambda_{U}(U_{0})=(U_{0})_{\rho}=\lambda_{U^{\prime}}(U_0);$$
then $$\lambda_{U^{\prime}}(x)\in (U_{0})_{\rho}\subseteq X_{0},$$
where there will be a contradiction. Hence, $\lambda_{U}$ and
$\lambda_{U^{\prime}}$ coincide on $U\cap U^{\prime}$ as mappings of
spaces.

$\emph{(c)}$ By gluing $\lambda_{U}$ along all such affine open
subsets $U$, we have a homeomorphism $\lambda$ of $X$ onto $X$ as a
topological space given by
\begin{equation*}
\lambda: x\in X \mapsto \lambda_{U}(x)\in X
\end{equation*}
where $x$ belongs to $U$ and $U$ is an affine open subset of $X$
such that $ \phi(U)$ is contained in some
 affine open subset $V$ of $Y$. That
is, $ \lambda|_{U}=\lambda_{U}. $ By $\emph{(b)}$ it is seen that
$\lambda $ is well-defined. It is clear that $\lambda$ is also an
automorphism of the scheme $(X,\mathcal{O}_{X})$.

Show that $\lambda$ is contained in $Aut\left( X/Y\right)$ such that
$ t\left(\lambda\right)=\rho$. In deed, as $\rho$ is an isomorphism
of $k(X)$ over $\phi^{\sharp}\left(k(Y)\right)$, it is seen that the
isomorphism $\lambda_{U}$ is over $Y$ by $\phi$ for any affine open
subset $U$ of $X$; then $\lambda$ is an automorphism of $X$ over $Y$
by $\phi$ such that $t\left(\lambda\right)=\rho$ holds.

This proves that there exists $\lambda\in Aut\left( X/Y\right) $
such that $t(\lambda)=\rho$ for each $\rho \in Gal\left( k\left(
X\right) /\phi^{\sharp}\left(k\left( Y\right)\right) \right) $.

Hence, ${t}$ is surjective. This completes the proof.
\end{proof}

\begin{corollary}
Let $X$ and $Y$ be arithmetic varieties and let $X$ be quasi-galois
closed over $Y$ by a surjective morphism $\phi$ of finite type.
Then there is a natural isomorphism
$$\mathcal{O}_{Y}\cong \phi_{\ast}(\mathcal{O}_{X})^{{Aut}\left( X/Y\right)}$$ where
$(\mathcal{O}_{X})^{{Aut}\left( X/Y\right)}(U)$ denotes the
invariant subring of $\mathcal{O}_{X}(U)$ under the natural action
of ${Aut}\left( X/Y\right)$ for any open subset $U$ of $X$.
\end{corollary}

\begin{proof}
Fixed any affine open sets $U_0$ of $X$ and $V_0$ of $Y$ with
$U_0\subseteq \phi^{-1}(V_0)$. By \emph{Theorem 3.26} we have
$$\phi^{\sharp}(\mathcal{O}_{Y}(V_0))= (\mathcal{O}_{X})^{{Aut}\left(
X/Y\right)}(U_0)= \phi_{\ast}(\mathcal{O}_{X})^{{Aut}\left(
X/Y\right)}(V_0)$$ since $k(X)=Fr(\mathcal{O}_{X}(U_0))$ and
$k(Y)=Fr(\mathcal{O}_{X}(V_0))$.

Now take any open set $V$ of $Y$ and put $U=\phi^{-1}(V)$. We must
have $$\phi^{\sharp}(\mathcal{O}_{Y}(V))=
(\mathcal{O}_{X})^{{Aut}\left( X/Y\right)}(U)=
\phi_{\ast}(\mathcal{O}_{X})^{{Aut}\left( X/Y\right)}(V).$$
Otherwise, if there is some element $w$ contained in the difference
set $(\mathcal{O}_{X})^{{Aut}\left( X/Y\right)}(U)\setminus
\phi^{\sharp}(\mathcal{O}_{Y}(V))$, we will have $$w\in
\mathcal{O}_{X}^{{Aut}\left( X/Y\right)}(U_1)\setminus
\phi^{\sharp}(\mathcal{O}_{Y}(V_1))$$ and then we will obtain a
contradiction, where $U_1\subseteq U$ and $V_1\subseteq V$ are
affine open sets such that $U_1\subseteq \phi^{-1}(V_1)$. This
completes the proof.
\end{proof}

\begin{remark}
Let $X$ and $Y$ be arithmetic varieties and let $X$ be quasi-galois
closed over $Y$ by a surjective morphism $\phi$ of finite type.
By \emph{Corollary 3.27} it is easily seen that the morphism $f$ must be affine.
\end{remark}

\subsection{Proof of the main theorem}

Now we are ready to prove the main theorem of the paper,
\emph{Theorem 2.1} in \S 2.

\begin{proof} (\textbf{Proof of Theorem 2.1.})
By \emph{Theorem 3.26} and \emph{Remark 3.28} it needs only to prove that $X$ is a pseudo-galois cover of $Y$ if $X$ and $Y$ have the same dimensions.

Let $\dim X= \dim Y$ and $G={Aut}\left( X/Y\right)$. It is clear
that $\phi$ is invariant under the natural action of ${Aut}\left(
X/Y\right)$ on $X$. By \emph{Corollary 3.27}, we have
$\mathcal{O}_{Y}\cong \phi_{\ast}(\mathcal{O}_{X})^{G}$. By
\emph{Theorem 3.26} again it is seen that $G$ is a finite group.
By \S 5 of \cite{VS1}, it is immediate
that $\phi$ is finite and then $X$ is a pseudo-galois cover of $Y$ by $\phi$.
\end{proof}

\begin{remark}
Fixed any  two arithmetic varieties $X$ and $Y$ such that  $X$ is quasi-galois closed over $Y$ by
 a surjective morphism $\phi$ of finite
type. By \emph{Theorem 2.1} it is seen that
the natural action of automorphism group ${Aut}\left( X/Y\right)$ on
the fiber $\phi^{-1}(y)$ is transitive at each $y\in Y$.

It follows that
 each point $x_0\in \phi ^{-1}(y_0)$ has a quasi-galois
set of geometric $y_{0}-$conjugations in $X$ for any point $y_0$ of
$Y$.
\end{remark}

\end{document}